\documentclass[12pt]{article}
\usepackage[letterpaper,top=2.5cm,bottom=2cm,left=3cm,right=3cm,marginparwidth=1.75cm]{geometry}
\usepackage[dvipsnames]{xcolor}
\usepackage{graphicx,amsfonts,amsmath,tikz-cd,amsthm}
\usepackage[final]{hyperref}
\usepackage{amssymb}
\usepackage{authblk}
\usepackage{tikz}
\usepackage{subfigure}
\usepackage[normalem]{ulem}
\usepackage{algpseudocode}
\usepackage[ruled,vlined,linesnumbered]{algorithm2e}

\usepackage{enumitem}

\title{Infinite Euclidean Distance Discriminants}
\date{}
\author[1]{Felix Rydell}
\author[2]{Emil Horobe\c{t}}
\affil[1]{ {\small Swedish Defence Research Agency, Stockholm, Sweden}}
\affil[2]{ {\small Sapientia Hungarian University of Transylvania, Romania}}

\newcommand{\new}[1]{{\color{black}{#1}}}

\providecommand{\keywords}[1]{\small \textbf{Keywords}: #1}
\providecommand{\subjclass}[1]{\small \textbf{Subject Class $2010$}: #1}

\newcommand{\rotsymbol}{\rotatebox{45}{\qedsymbol}}

\usepackage{listings}
 
\definecolor{codedarkgreen}{RGB}{51, 133, 4}
\definecolor{codemaroon}{RGB}{133, 5, 63}
\definecolor{codeteal}{RGB}{0, 128, 96}
\lstdefinelanguage{Macaulay2}{
basicstyle=\small\ttfamily,
alsoletter=",
classoffset=1,
keywords={matrix,minors,gb,transpose,det,ideal,apply,subsets,ker,gens,fold,flatten,entries},
keywordstyle={\color{blue}},
classoffset=2,
morekeywords={from, to, list},
keywordstyle={\color{codemaroon}},
classoffset=3,
morekeywords={QQ},
keywordstyle={\color{codedarkgreen}},
classoffset=4,
morekeywords={MonomialOrder},
keywordstyle={\color{codeteal}},
xleftmargin=1.5cm,
xrightmargin=1em,
columns=fullflexible,
keepspaces=true,
stepnumber=1,
numbers=none,
captionpos=b,
showspaces=false,
frame=none
}

\usepackage[capitalise, noabbrev, nameinlink]{cleveref}

\crefname{listing}{Code}{Codes}
\Crefname{listing}{Code}{Codes}
\crefname{lstlisting}{Code}{Codes}
\Crefname{lstlisting}{Code}{Codes}

\theoremstyle{plain}
\newtheorem{theorem}{Theorem}[section] 
\newtheorem{proposition}[theorem]{Proposition}
\newtheorem{lemma}[theorem]{Lemma}

\newtheorem{corollary}[theorem]{Corollary}

\theoremstyle{definition}
\newtheorem{definition}[theorem]{Definition} 
\newtheorem{example}[theorem]{Example}
\newtheorem{remark}[theorem]{Remark}

\theoremstyle{plain}
\newtheorem{introtheorem}{Theorem} 
\newtheorem{introproposition}[introtheorem]{Proposition}

\newcommand{\PP}{\mathbb{P}}
\newcommand{\RR}{\mathbb{R}}

\newcommand{\CC}{\mathbb{C}}

\newcommand{\QQ}{\mathbb{Q}}

\begin{document}

\maketitle

\begin{abstract}
 We study infinite Euclidean distance discriminants of algebraic varieties, defined as the loci of data points whose fibers under the second projection from the Euclidean distance correspondence are positive-dimensional. In particular, these discriminants contain all data points with infinitely many critical points for the nearest-point problem. We present computer code that computes the infinite Euclidean distance discriminant, and use it to present numerous varieties with nonempty such discriminants. Moreover, we prove that for any data point, the fiber under the second projection is contained in a finite union of hyperspheres centered at that point. For curves, we include a complete characterization; their infinite Euclidean distance discriminants turn out to be affine linear spaces. Finally, we introduce and characterize skew-tube surfaces in three-dimensional space. By construction, these have a one-dimensional infinite Euclidean distance discriminant. We further demonstrate that many skew-tubes have significantly lower Euclidean distance degrees than generic surfaces of the same degree.
\end{abstract}
. \keywords{Parametric optimization, Euclidean Distance Degree, Constrained Critical Points, Euclidean Distance Discriminant, Purity of branch locus}

%\tableofcontents

\section*{Introduction} 
In many real-world applications, minimizing the distance between a data point \new{derived from a noisy measurement} and a target model plays a crucial role. Such problems are known as \textit{Euclidean distance problems} or \textit{nearest point problems}, and appear, for example, in computer vision~\cite{hartley1997triangulation}, geometric modeling~\cite{thomassen2004closest}, and low-rank matrix approximation~\cite{chu2003structured}. From an algebraic point of view, they form a main topic of metric algebraic geometry~\cite{breiding2024metric}. 

If the model is an algebraic variety, Draisma et al.~\cite{draisma2016euclidean} showed that, for a generic data point, the number of complex critical points of the distance function is constant. This number is the Euclidean distance degree (ED degree) of the variety. In concrete applications, however, one is typically interested only in the real solutions. The number of real solutions is not constant for generic data. Instead, there exists a variety, called the Euclidean distance discriminant (ED discriminant)---classically the focal locus~\cite{catanese2000focal}, or the \textit{evolute} for plane curves---which partitions real space into connected regions, where the number of real critical points is constant. This fact can be used in implementations of homotopy continuation to solve the critical equations by only tracking real solutions rather than all complex ones.

\begin{figure}
    \centering
   \includegraphics[width=0.45\textwidth]{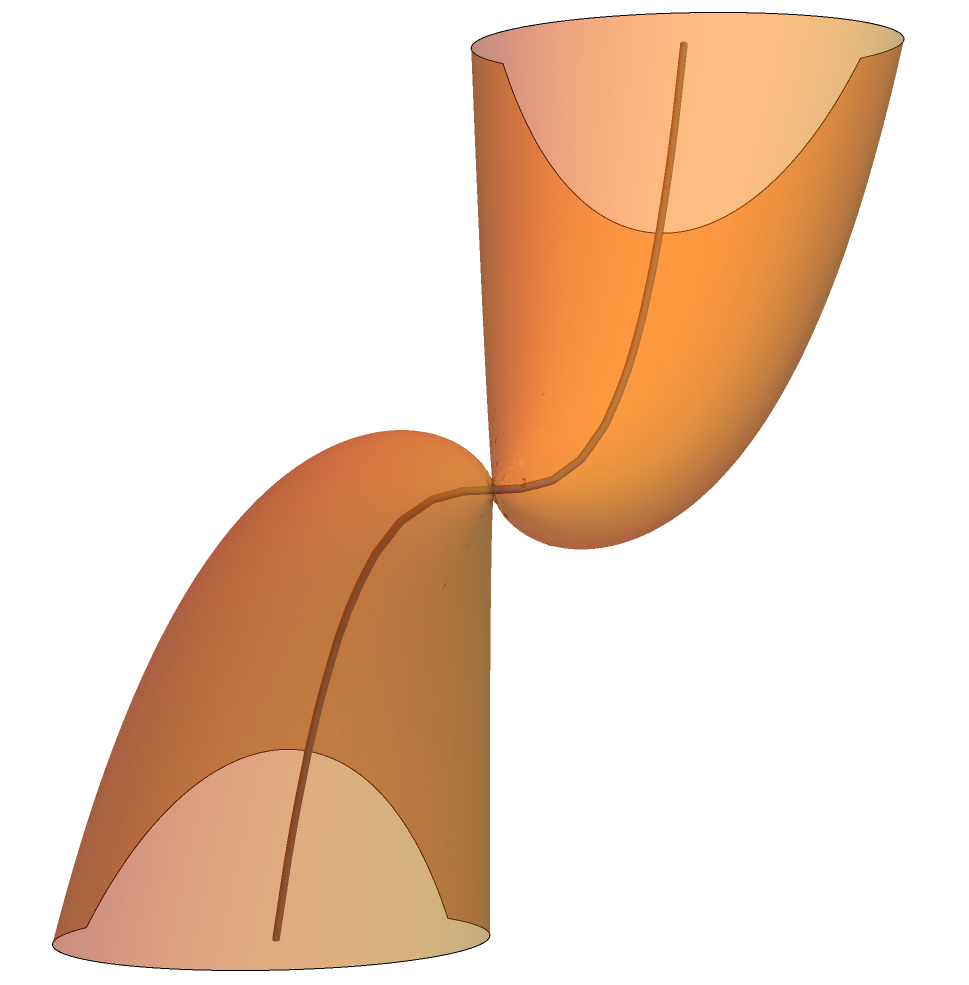} \hspace{1em}\includegraphics[width=0.45\textwidth]{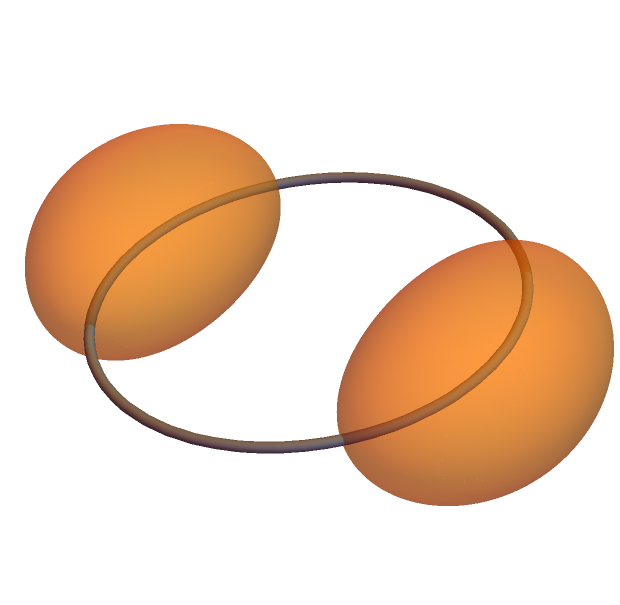}
    \caption{Surface\new{s (orange)} whose infinite ED discriminants contain curves (gray). On the left: the degree-10 surface (T1) and a twisted cubic (\Cref{ex: twisted skewtubes}). On the right: the degree-8 surface (C1) and a circle (\Cref{ex: circle skewtubes}).}
    \label{fig: teaser}
\end{figure}

Given an affine variety $X$ and a fixed data point $u$, the number of distinct complex critical points may differ from the ED degree in several ways:
\begin{enumerate}
    \item Two critical points may collide, so that a critical point acquires multiplicity of at least two. 
    \item The number of critical points may exceed the ED degree, becoming infinite. 
    \item Some critical points may lie in the singular locus (the set of such $u$ is the singular data locus~\cite{horobect2017data}) or at infinity (the set of such $u$ is the atypical discriminant~\cite{JST24}). 
\end{enumerate}
\new{In this paper, we study the second scenario in detail. Recall} that the Euclidean distance correspondence (ED correspondence) is the closure of the pairs $(x,u)$, where $x$ \new{is critical} to $u$. We define the \textit{infinite Euclidean distance discriminant} (infinite ED discriminant) of $X$, denoted $\Sigma^{\infty}_X$, as the locus of data points $u$ whose fiber $\mathcal E_X(u)$ under the second projection $(x,u)\mapsto u$ is positive-dimensional. Illustrations of varieties with nonempty infinite ED discriminants are shown in \Cref{fig: teaser,fig: examples,fig: skewtubes}, and were created using \texttt{Mathematica} \cite{Mathematica}.

The structure and main results of this paper are as follows. In \Cref{s: prelim}, we introduce notation, definitions, and well-established results that we use throughout. In \Cref{s: Comp}, we present pseudo-code (\Cref{alg: sigma}) to compute infinite ED discriminants. We provide accompanying \texttt{Macaulay2} \cite{M2} in the \new{source folder of the} \texttt{arXiv} version of the paper. We show the following: 
\begin{introtheorem}[Code for infinite ED discriminants] \Cref{alg: sigma} terminates and correctly returns the infinite Euclidean distance discriminant.
\end{introtheorem}
In \Cref{s: examples}, we implement \Cref{alg: sigma} to efficiently compute infinite ED discriminants of several of varieties. We provide families of varieties with nonempty infinite ED discriminants, such as varieties contained in a hypersphere, offset hypersurfaces, and varieties of revolution (\Cref{thm: rev and tubes}). In \Cref{s: Properties}, we study general properties of infinite ED discriminants. Our main theorem is: 

\begin{introtheorem}[Structure of critical points]\label{thm: intro spheres} Let $X\subseteq \CC^n$ be an irreducible variety. For any $u\in \CC^n$, $\mathcal E_X(u)$ is contained in a finite union of hyperspheres centered at $u$. 
\end{introtheorem}

In \Cref{s: curves}, we classify infinite ED discriminants for curves. As a consequence, a generic complete intersection curve has an empty infinite ED discriminant. We believe that this phenomenon extends to generic complete intersections of arbitrary dimension, although we do not provide a proof. 

\begin{introtheorem}[Classification for curves]\label{thm: intro curves} Let $X$ be an irreducible curve in $\CC^n$. 
\begin{enumerate}[label=(\Roman*)]
    \item If $X$ is not contained in a hypersphere, then $\Sigma_X^\infty=\emptyset$.
    \item Otherwise, let $c$ be the center of such a hypersphere. If $V$ is the smallest affine linear space containing $X$, and $N\:V$ is its normal space, then 
    \begin{align}
     \Sigma_X^\infty=\{c+y: y\in N\: V\}.  
    \end{align}
\end{enumerate}
\end{introtheorem}

In $\CC^3$, we say that a surface $X$ is a \textit{skew-tube} around a curve $C$ if, for generic $c\in C$, the set of critical points $C_X(c)$ is one-dimensional, and the union of these is dense in $X$ (see \Cref{def: skew} for the general definition). These skew-tubes are characterized by associating one (or possibly multiple) squared radii $R$ to each point $c$ of $C$. Algebraically, this relation is encoded by a polynomial equation $f(c,R)=0$, which defines a curve $\mathcal U\subseteq \CC^3\times \CC$ in the coordinates $(c,R)$. Letting $\mathrm{pr}_1$ denote the first projection $(c,R)\mapsto c$, we obtain $C=\overline{\mathrm{pr}_1(\mathcal U)}$. \new{Locally, $\mathcal U$ is cut out by 3 polynomials $f_1,f_2,f_3$ whose gradients span the normal space of $\mathcal U$ for generic $(c, R)$. We obtain a \textit{local tangent} $T_{(c,R)}^{\mathrm{loc}} \:\mathcal U$  by applying Cramer's rule on the Jacobian of $f_i$. To proceed, let $x\cdot y=x_1y_1+\cdots +x_ny_n$ denote the standard Euclidean complex inner product.}

\begin{introtheorem}[Construction of skew-tubes]\label{thm: intro construction} Let $\mathcal U\subseteq \CC^3\times \CC$ be an irreducible curve that is cut out locally by $f_1,f_2,f_3\in \CC[c_1,c_2,c_3,R]$\new{, which} determine a local tangent $T_{(c,R)}^{\mathrm{loc}}\:\mathcal U$. Assume that $C=\overline{\mathrm{pr}_1(\mathcal U)}$ is a curve. Define the ideal $I^{\mathcal U}\subseteq S:=\CC[x_1,x_2,x_3,c_1,c_2,c_3,R]$ as the sum of the defining ideal of $\mathcal U$ and the ideal generated by the conditions
\begin{align}
   (x-c)\cdot (x-c)=R, \quad T_{(c,R)}^{\mathrm{loc}}\:\mathcal U\cdot \left(x-c,\frac{1}{2}\right)=0.
\end{align}
Let $I_{\mathrm{sat}}^{\mathcal U}$ denote the saturation of $I^{\mathcal U}$ with respect to the ideal generated by $T_{(c,R)}^{\mathrm{loc}}\:\mathcal U=0$. Then any component $X$ of the zero locus of
\begin{align}
   I_{\mathrm{elim}}^{\mathcal U}:=I_{\mathrm{sat}}^{\mathcal U}\cap \CC[x_1,x_2,x_3]
\end{align}
satisfies $C\subseteq\Sigma_X^\infty$. 

Moreover, for any surface $X$ in $\CC^3$ that is a skew-tube around a curve $C$, there is an irreducible curve $\mathcal U\subseteq \CC^4$ with $C=\overline{\mathrm{pr}_1(\mathcal U)}$ such that $X$ is a component of the zero locus of $I_{\mathrm{elim}}^{\mathcal U}$.
\end{introtheorem}

We demonstrate that skew-tubes exhibit desirable optimization properties. In particular, their critical ideals can be simplified, enabling the efficient computation of critical points even for high-degree surfaces.

\begin{introproposition}[Optimization for skew-tubes]\label{introprop:opt} Assume that the surface $X\subseteq \CC^3$ is a skew-tube around the curve $C$. Let $\mathcal U\subseteq \CC^3\times \CC$ be as in \Cref{thm: intro construction}. If $\mathcal U$ is locally cut out by polynomials $f_1,f_2,f_3$ of degrees $d_1, d_2, d_3$, then  
\begin{align}
    \mathrm{EDdeg}(X) \le d_1d_2d_3\big(2(d_1+d_2+d_3)-3\big).
\end{align}
\end{introproposition}
We compare this bound with the expected ED degree (\cite[Proposition~2.6]{draisma2016euclidean}) in several examples and observe a substantial discrepancy. For instance, we construct a degree-$10$ skew-tube (denoted (T2) in \Cref{ex: twisted skewtubes}) whose ED degree is bounded above by 28 (\Cref{introprop:opt}), whereas the expected ED degree for a generic surface of degree $10$ is \new{$910$}. By exploiting the structure of skew-tubes, we obtain a simplified critical ideal \new{for (T2)} that enables us to compute the critical points of a data point in under a hundredth of a second, compared to more than 45 minutes using the standard critical equations~\cite[Section 2]{draisma2016euclidean}. These computations indicate that the actual ED degree \new{of (T2)} is $12$ (\Cref{ex: skew ED}). Together with the other examples, this \new{suggests} a correlation between a nonempty infinite ED discriminant and a low ED degree.

It is pertinent to place our results in context. A weaker version of \Cref{thm: intro spheres} follows from a result in the PhD thesis of Sodomaco~\cite[Corollary~4.3.7]{sodomaco2020distance}, where $X$ is assumed to be transerval to the isotropic quadric. Furthermore, in work of Catanese and Trifogli, \Cref{thm: intro spheres} appears for smooth hypersurfaces~\cite[Section~6, Proposition~2]{catanese2000focal}, and \Cref{thm: intro construction} can be viewed as a special case of~\cite[Section~6, Definition~6]{catanese2000focal}. 
In contrast to the work of Catanese and Trifogli, our perspective on this topic is from an applied and computational point of view. In recent work on the atypical discriminant, Joi{\c{t}}a, Siersma, and Tib{\u{a}}r proved \Cref{thm: intro curves} in the case $n=2$~\cite[Theorem 2.5 (c)]{JST24}. Scenarios with infinitely many critical points have been considered by Beorchia, Gallet, and Logar~\cite{BGL24} for eigenvectors of ternary cubics, and by Abo, Portakal, and Sodomaco in the study of Nash equilibria~\cite{APS25}.

\bigskip

\paragraph{\textbf{Acknowledgements}.} Felix Rydell was supported by the Knut and Alice Wallenberg Foundation within their WASP (Wallenberg AI, Autonomous
Systems and Software Program) AI/Math initiative. Emil Horobe\c{t} was supported by the project “Singularities and Applications” - CF 132/31.07.2023 funded by the European Union - NextGenerationEU - through Romania’s National Recovery and Resilience Plan. The authors are grateful to Mike Stillman for conveying the main ideas for \Cref{alg: sigma}, to Luca Sodomaco for providing \Cref{ex: twisted}, and to Kathlén Kohn for her input and feedback. \new{We additionally thank Sofía Errazuriz and Thomas Georg Grill for their input on the proofs of \Cref{prop: bisec} and \Cref{thm: deg EDpoly}.}

%%%%%%%%%%%%%%%%%%%%%%%%%%%%%%%%%%%%%%%%%%%%%%%%%%%%%%%%%%%%%%%%%%%%%%%%%%%%%%%%%%%%%%%%%%%%%%%%%%%%

\section{Preliminaries}\label{s: prelim}
 \setcounter{equation}{0}
\numberwithin{equation}{section} 
Let $X_{\mathbb{R}}\subseteq \mathbb{R}^n$ be a real irreducible variety, and take a data point $u\in\mathbb{R}^n$ in the ambient space. Consider the following constrained optimization problem:
\begin{equation}\label{MinProblem}
\begin{cases}
\text{minimize } d_u(x):=\displaystyle\sum_{i=1}^n (u_i-x_i)^2, \\
\text{subject to }x\in X_\RR.
\end{cases}
\end{equation}
Here, $d_u$ denotes the squared Euclidean distance function. To apply techniques from \new{algebraic geometry}, we pass to the complexification $X\subseteq \mathbb{C}^n$, defined by the same real polynomials as $X_{\mathbb{R}}$. Over the complex numbers, the minimization is not well-defined, but we can nevertheless study the associated critical equations and their zero locus, that is, the set of critical points.

Suppose that the defining ideal of $X$ is generated by the polynomials $f_1,\ldots,f_s$. A classical approach to solving constrained optimization problems like \eqref{MinProblem} is to use Lagrange multipliers~\cite{boyd2004convex}. \new{This method uses the fact that a smooth optimum $x$ satisfies}
\begin{equation}\label{eq: CritEqu}
\begin{cases}
\displaystyle\nabla d_u(x)+\sum_{i=1}^s \lambda_i \nabla f_i(x)=0,\\
f_1(x)=\cdots =f_s(x)=0,
\end{cases}
\end{equation}
for some \new{multipliers} $\lambda_i\in \CC$. At a smooth point $x\in X$, the normal space $N_x\: X\subseteq \CC^n$ is the vector space spanned by the gradients \new{of $f_i$}, and the tangent space $T_x\: X=\{v\in \CC^n: v \cdot  w=0 \textnormal{ for every }w\in N_x\: X\}$. Here and throughout this paper, $v\cdot w = v_1w_1+\cdots v_nw_n$ is the standard complex dot product. For singular $x\in X$, the gradients span a space of lower than expected dimension, which changes the character of the optimization problem \eqref{eq: CritEqu}. Therefore, we restrict to the locus of smooth \new{points}, denoted $X_{\mathrm{reg}}$, and exclude the locus of singular points, denoted $X_{\mathrm{sing}}$. Using this notation, \eqref{eq: CritEqu} is on $X_{\mathrm{reg}}$ equivalent to
\begin{align}\label{eq: reg cond}
    \begin{cases}
u-x \in N_x\:X,\\
x\in X_{\mathrm{reg}}.
\end{cases}
\end{align}
\new{Given a data point $u$, the set of $x\in X$ satisfying  \eqref{eq: reg cond} are the \textit{critical points} to $u$.} The closure $\mathcal E_X$ of all pairs $(x,u)$ that satisfy \eqref{eq: reg cond} is the \textit{ED correspondence}:
\begin{align}\begin{aligned}
   \mathcal E_X:=\overline{E_X}\quad \textnormal{ where }\quad  E_X:=\{(x,u)\in \CC^n\times \CC^n: \eqref{eq: reg cond} \textnormal{ holds}\}.
\end{aligned}
\end{align}
The variety $\mathcal{E}_X$ is irreducible of dimension $n$ \cite[Theorem 4.1]{draisma2016euclidean}. There are two natural projections $\mathrm{pr}_i$ from $\mathcal E_X$ onto its factors:
\begin{equation}\label{eq: EX proj}
\begin{tikzcd}
&\mathcal{E}_X \arrow{ld}[swap]{\mathrm{pr}_{1}}\arrow{rd}{\mathrm{pr}_2}& \\
X & & \mathbb{C}^n
\end{tikzcd}
\end{equation}
Assuming that the second projection is dominant, it is generically finite-to-one, and the cardinality of the fiber over a generic data point $u\in \CC^n$ is, by definition, the \textit{ED degree} of $X$, written $\mathrm{EDdeg}(X)$. For affine cones, the \textit{ED discriminant} is the branch locus of $\mathrm{pr}_2:\mathcal P\mathcal E_X\to \CC^n$, where $ \mathcal P\mathcal E_X$ is the projective ED correspondence~\cite[Section 7]{draisma2016euclidean}. It is stated in the literature~\cite[Section 7]{draisma2016euclidean}, \cite[Chapter 1, Section 2]{sodomaco2020distance} that the ED discriminant is the closure of the set of data points with less than ED degree many distinct critical points. %https://stacks.math.columbia.edu/tag/02V5?utm_source=chatgpt.com % https://stacks.math.columbia.edu/tag/02NI?utm_source=chatgpt.com
\new{However, the differential of $\mathrm{pr}_2$ is also rank-deficient when the fiber is positive-dimensional.} Accordingly, we define the (\textit{affine}) ED discriminant as 
\begin{align}\label{eq: SigmaX}
\Sigma_X := 
\left\{\, u \in \CC^n \;:\;
\begin{aligned}
    & \textnormal{there is an isolated solution } (x,u) \in \mathcal E_X \\
    &\textnormal{of multiplicity at least 2 in } x, \\
    &\textnormal{or there are infinitely many solutions } (x,u) \in \mathcal E_X
\end{aligned}
\right\}.
\end{align}
In this paper, we study the subvariety of $\Sigma_X$ of points with positive-dimensional fibers:

\begin{definition}\label{def: infinite ED disc}
Let $X$ be an irreducible variety. The \textit{infinite Euclidean distance discriminant} (infinite ED discriminant) of $X$ is 
\begin{align}
    \Sigma_X^{\infty}:=\{u\in \mathbb{C}^n: \dim \mathrm{pr}_2^{-1}(u)\ge 1\}.
\end{align}
We denote $\mathrm{pr}_2^{-1}(u)$ by $\mathcal E_X(u)$, and we denote $C_X(u)$ the set of critical points, meaning all $x$ such that $(x,u)$ satisfy \eqref{eq: reg cond}. 
\end{definition}

By definition, $C_X(u)\subseteq \mathcal E_X(u) $, and as a consequence, $\overline{C_X(u)}\subseteq \mathcal E_X(u)$. A natural question is whether this inclusion is, in fact, always an equality. In \Cref{ex: volcano}, we show an example where this is not the case.

\begin{remark} The upper semi-continuity of dimensions \cite[Chapter 1, \S 8, Corollary 3]{mumford1999red}) ensures that $\Sigma_X^{\infty}$ is a variety.
\end{remark}

\begin{lemma} Let $X\subseteq \CC^n$ be irreducible. If $\mathrm{EDdeg}(X)>0$, then $\dim \Sigma_X^\infty \le n-2$.
\end{lemma}

\begin{proof} The preimage $\mathrm{pr}_2^{-1}(\Sigma_X^\infty)$ is of greater dimension than $\Sigma_X^\infty$ by construction, and it lies in the ramification locus of \new{$\mathrm{pr}_2$}. Since the ramification locus is at most a hypersurface, the dimension of $\mathrm{pr}_2^{-1}(\Sigma_X^\infty)$ is at most $n-1$. 
\end{proof}

%%%%%%%%%%%%%%%%%%%%%%%%%%%%%%%%%%%%%%%%%%%%%%%%%%%%%%%%%%%%%%%%%

The closure of the image of 
\begin{align}
    \mathcal E_X\to \CC^n\times \CC, \quad (x,u)\mapsto (u,d_u(x)),
\end{align}
is irreducible, and of dimension $n$ if and only if  $\mathrm{EDdeg}(X)>0$, which we assume from now on. \new{We call this variety the \textit{offset hypersurface} and denote it by $\mathcal O_X$.} It is cut out by a single irreducible polynomial $F_u(R)\in \CC[u,R]$. The \textit{ED polynomial} of $X$ is $f_u(\epsilon):=F_u(\epsilon^2)\in \CC[u,\epsilon]$~\cite{horobect2019offset}. It records information such as the (complex) distances $\epsilon$ between generic $u$ and its critical points, the ED degree, and the second fundamental form~\cite[Theorem 7.21]{breiding2024metric}. For a fixed $\epsilon$, we also define the \textit{$\epsilon$-offset hypersurface} $\mathcal O_{X,\epsilon}$ as 
\begin{align}
\mathcal O_{X,\epsilon}:=\{u\in \CC^n: f_u(\epsilon)=0\}.
\end{align}
\new{One can check that $\mathcal O_{X,\epsilon}$ is indeed a hypersurface for any $\epsilon$, meaning it is neither the ambient space or empty, under the assumption that the ED degree is positive.}

The \textit{bisector locus} of a variety is the closure of the set of data points $u$ that admit distinct critical points $x,y$ equidistant to $u$. Formally,
\begin{align}
B_X := 
\overline{\left\{\, u \in \CC^n \;:\;
\begin{aligned}
    &\textnormal{there exist distinct }x,y\in X\textnormal{ such that }\\
    &(x,u),(y,u)\in \mathcal E_X \textnormal{ and }d_u(x)=d_u(y)
\end{aligned}
\right\}}.
\end{align}
By definition, the infinite ED discriminant is contained in both the ED discriminant and the bisector locus:
  \begin{align}
    \Sigma_{X}^{\infty}\subseteq \Sigma_X\cap B_{X}.
    \end{align}

\begin{proposition}\label{prop: bisec} Let $X\subseteq \CC^n$ be an irreducible variety. The bisector locus is at most a hypersurface.  
\end{proposition}
We postpone the proof to the next subsection.

\begin{remark} It has been reported~\cite{horobect2019offset,breiding2024metric} that the \textit{offset discriminant} $\Delta_X$\new{, defined as the $\epsilon^2$-discriminant of the ED polynomial,} is the union of the bisector locus and the ED discriminant. This fails in the affine setting: the offset discriminant contains a data point $u$ if and only if there exists a sequence of generic data points $u_n\to u$ (in Euclidean topology) such that one of the following holds:
\begin{enumerate}
\item For some $\epsilon \in \CC$, there are sequences $x_{i,n}\in \mathcal E_X(u_n)$, for $i=1,2$, with $x_{1,n}\neq x_{2,n}$ for each $n$ and 
\begin{align}
  \new{d_{u_n}(x_{i,n})}\to \epsilon^2 \quad \textnormal{for }i=1,2. 
\end{align}
\item For each $\epsilon\in \CC$, there is a sequence $x_n\in \mathcal E_X(u_n)$ with  
\begin{align}
    \new{d_{u_n}(x_n)}\to \epsilon^2.
\end{align}
\end{enumerate}
In particular, if $\mathcal E_X(u)$ has an isolated double point at infinity, it is not guaranteed a priori to be captured by the bisector locus or the ED discriminant. See for instance~\Cref{ex: foot}.
\end{remark}

%%%%%%%%%%%%%%%%%%%%%%%%%%%%%%%%%%%%%%%%%%%%%%%%%%%%%%%%%%%%%%%%%

\subsection{Local parametrizations}\label{ss: local par} In this paper, we use local smooth \new{parametrizations} of algebraic varieties. Such \new{parametrizations} need not be given by polynomials. This is the case for the \new{nonrational} Fermat curve
\begin{align}
    x^3+y^3+1=0.
\end{align}
%https://en.wikipedia.org/wiki/Rational_variety
The smooth locus of any variety is, however, a manifold~\cite[Chapter I, Theorem 5.3]{hartshorne2013algebraic}. This means that there is a local smooth parametrization around any smooth point of a variety. Our first application of smooth local parametrization is the next lemma.

\begin{lemma}\label{le: offset duality} Let $u\in \mathcal O_{X,\epsilon}$ be smooth and let $u(s)$ be a local smooth parametrization \new{around $u(0)=u$}. Let $x(s)\in X$ be a local smooth function such that \new{$d_{u(s)}(x(s))=\epsilon^2$} and $x(s)$ is critical to $u(s)$, for each $s$. Then
\begin{align}
u(s)-x(s)\in N_{u(s)}\:\mathcal O_{X,\epsilon},
\end{align}
meaning that $u(s)\in \mathcal O_{X,\epsilon}$ is critical to $x(s)$ (viewed as a data point).
\end{lemma}

\begin{proof} Deriving both sides of \new{$d_{u(s)}(x(s))=\epsilon^2$} with respect to $s_i$, we have
\begin{align}
    \Big(\frac{\partial}{\partial s_i} u(s)- \frac{\partial}{\partial s_i} x(s)\Big)\cdot (u(s)-x(s))=0.
\end{align}
The fact that $u(s)-x(s)\in N_{x(s)}\: X$ means exactly that $\frac{\partial}{\partial s_i} x(s) \cdot (u(s)-x(s))=0$ for each $s_i$. Therefore, we are left with $ \frac{\partial}{\partial s_i} u(s)\cdot (u(s)-x(s))=0 $ for each $s_i$, meaning that $u(s)-x(s)\in N_{u(s)}\: \mathcal O_{X,\epsilon}$.
\end{proof}

Let $I\subseteq \CC[x,u]$ be a prime ideal and assume that $I\cap \CC[u]=\langle 0\rangle$, i.e., the projection $(x,u)\mapsto u$ from the zero locus $V(I)\subseteq \CC^n\times \CC^n$ is dominant \new{onto $\CC^n$}. It is a fact of algebraic geometry that for generic $p\in \CC^n$, the \textit{specialization ideal} $I_p:=\{f(x,p): f\in I\}\subseteq \CC[x]$ is radical. In particular, if, for generic $p$, $V(I_p)$ is zero-dimensional, then it consists of distinct points of multiplicity 1. In this case, the implicit function theorem applies: for $x\in V(I_p)$ and a smooth local parametrization $p(s)\in \CC^n$ of $p(0)=p$, there is a unique smooth function $x(s)$ such that $x(s)\in V(I_{p(s)})$ for every small $s$. The point $x\in V(I_p)$ has multiplicity one if it is smooth in $V(I_p)$, which is equivalent to
\begin{align}
\mathrm{rank}\:\mathrm{Jac}_x(I_p) = n,
\end{align}
where $\mathrm{Jac}_x$ denotes the Jacobian with respect to the $x$-variables. The cardinality of $V(I_p)$ is constant for generic $p\in \CC^n$~\cite{borovik2025short}. If the zeros of $I_p$ (counting multiplicities) exceed this number, then $V(I_p)$ is positive-dimensional.

\begin{proof}[Proof of \Cref{prop: bisec}] If $\mathrm{EDdeg}(X)$ is $0$ or $1$, then there is nothing to prove; for generic $u\in \CC^n$, there is at most one critical point. If $\mathrm{EDdeg}(X)\ge 2$, then we assume by contradiction that the bisector locus equals $\CC^n$. As $\CC^n$ has ED degree 1, $X$ is a proper subvariety of $\CC^n$. Consider the following disjoint union that partitions $\mathcal E_X$:
\begin{align}
    \mathcal E_X = \bigcup_{\epsilon\in \CC}\mathcal E_{X,\epsilon} \quad \textnormal{where} \quad \mathcal E_{X,\epsilon}:= \mathcal E_X \cap \{(x,u): \new{d_u(x)}=\epsilon^2\}.
\end{align}
There are two cases:
\begin{enumerate}
    \item $\mathcal E_X$ equals $\mathcal E_{X,\epsilon}$ for some $\epsilon$. Since $(x,x)\in \mathcal E_X$ for smooth $x\in X$, this $\epsilon$ must equal $0$.\label{enum: EX equals}
    \item \new{$\mathcal E_{X,\epsilon}$ are} hypersurfaces in $\mathcal E_X$ for generic $\epsilon$.\label{enum: all hyper}
\end{enumerate}

\new{Since $B_X$ is assumed to equal $\CC^n$, the implicit function theorem (and the discussion before the proof) implies the following. There are, for generic $u\in \mathcal \CC^n$, distinct local smooth functions $x_1(u), x_2(u)$ in $X$ of critical points to the data $u$ such that 
\begin{align}\label{eq: two crit}
    d_u(x_i(u))=\epsilon(u)^2 \quad \textnormal{ and }\quad u-x_i(u)\in N_{x_i(u)}\:  X,
\end{align}
for $i=1,2$, and some smooth $\epsilon : \CC\to \CC$. In the case of \Cref{enum: EX equals}, $\epsilon(u)=0$, and an argument similar to the proof of \Cref{le: offset duality} shows that $u-x_i(u)\in N_{u}\: \mathcal O_{X,0}$. However, by assumption, $\mathcal O_{X,0}=\CC^n$, which implies the contradiction $x_1(u)=x_2(u)$.}

\new{In the case of \Cref{enum: all hyper}, $\epsilon(u)$ is generic in $\CC$ for generic $u\in \CC^n$ (otherwise $F_u(R)$ must by irreducibility equal $R-\epsilon^2$ up to scaling, for some $\epsilon$, which is \Cref{enum: EX equals}). Now fix a generic $u\in \CC^n$. The associated $x_1(u)$ from \eqref{eq: two crit} is generic in $X$. In particular, $(x_1(u), u)$ is a smooth point of $\mathcal E_{X,\epsilon}$, where $\epsilon:=\epsilon(0)$. Let $(x(t), u(t)), t\in \CC^{n-1},$ be a local parametrization of $\mathcal E_{X,\epsilon}$ around $(x(0),u(0))=(x_1(u),u)$. Then $x(t)$ must coincide with $x_1(t):=x_1(u(t))$. Denoting $x_2(t):=x_2(u(t))$, we have that $d_{u(t)}(x_i(t))=\epsilon^2$, for $i=1,2$, and every $t$ close to $0$. By \Cref{le: offset duality}, we find that $u(t)-x_i(t)\in N_{u(t)}\: \mathcal O_{X,\epsilon}$, for $i=1,2$. Note that $\mathcal O_{X,\epsilon}$ is a hypersurface, and therefore $N_{u(t)}\: \mathcal O_{X,\epsilon}$ is spanned by $u(t)-x_i(t)$ for both $i$. 

Consider the normal line $N(u,\epsilon)=\{u+ v: v\in N_{u}\: \mathcal O_{X,\epsilon}\}$. It contains the critical points $x_1:=x_1(0)$ and $x_2:=x_2(0)$. Note that there are precisely two (counting multiplicities) solutions to the quadratic equation $d_u(z)=\epsilon^2$ in $z\in N(u,\epsilon)$, since $\epsilon$ is generic and nonzero. These correspond to $x_1$ and $x_2$, and we deduce the relation $x_2=2u-x_1$. 

To reach a contradiction, consider data points $\hat u(\lambda):=u + (\lambda-1)(u-x_1)$, for $\lambda\in \CC$, and note that $d_{\hat u(\lambda)}(x_1)=(\lambda\epsilon)^2$. The point $x_1\in X$ is critical to $\hat u(\lambda)$ and $\hat u(\lambda) \in N(u,\epsilon)$, since $ N_{u}\: \mathcal O_{X,\epsilon}$ is the span of $u-x_1$. For generic $\lambda$, we can repeat the above to find that $\hat u(\lambda) -x_1$ spans $N_{\hat u(\lambda)}\: \mathcal O_{X,\lambda\epsilon}$, and as a consequence, $N_{\hat u(\lambda)}\: \mathcal O_{X,\lambda\epsilon} = N_{u}\: \mathcal O_{X,\epsilon}$. In particular, we have $N(\hat u(\lambda), \lambda \epsilon) = N(u, \epsilon)$. For $x_2(u)$ as in \eqref{eq: two crit}, define $\hat{x}_2(\lambda):= x_2(\hat{u}(\lambda))$. Then $\hat{x}_2(\lambda)$ is critical to $\hat{u}(\lambda)$ at squared distance $(\lambda\epsilon )^2$ and satisfies $\hat x_2(\lambda)=2\hat u (\lambda)-x_1$. The map $\hat x_2:\CC \to N(u,\epsilon)$ is dominant, implying that the whole line $N(u,\epsilon)$ is contained in $X$. In particular $u\in X$, which is a contradiction, since $X\neq \CC^n$ and $u$ was generic.} 
\end{proof}

We use \new{parametrizations} with additional structure. For this, we refer to Ehremann's fibration theorem~\cite[Section 8.5]{dundas2018short}:

\begin{theorem} Let $f: E \to  M$ be a
proper submersion. Then $f$ is a locally trivial fibration.
\end{theorem}
The map $f$ is \textit{proper} if the inverse image of every compact set is compact. It is a \textit{submersion} if it is smooth, and the differential map is of full rank. We consider the Euclidean topology when we use this theorem: $f$ is a \textit{locally trivial fibration} if for each $p\in M$, there is a Euclidean open neighborhood $U$ around $p$ and a diffeomorphism $h$ making the following diagram commute: 
\begin{equation}\label{eq: Ehresmann proj}
\begin{tikzcd}
f^{-1}(U)\arrow[rr, "h"] \arrow[rd,  "f\restriction_{f^{-1}(U)}"'{near start}]& & U\times f^{-1}(p) \arrow[ld, "\mathrm{pr}_U"{pos=0.4, xshift=3pt, yshift=-1.5pt}]\\
 & U & 
\end{tikzcd}
\end{equation}

%https://en.wikipedia.org/wiki/Rank_(differential_topology)

\subsection{Notes on smoothness} 
A natural way to partition a variety $X$ is via the chain of inclusions
\begin{align}
    X\supseteq X_{\mathrm{sing}}\supseteq (X_{\mathrm{sing}})_{\mathrm{sing}}\supseteq \cdots.
\end{align}
However, not all points in the differences between consecutive varieties in this chain necessarily look locally the same. A classic example is the Whitney umbrella. A more fine-grained version of \new{the} partition is the Whitney stratification~\cite{trotman2020stratification}, which we recall from~\cite[Section 4.2]{nicolaescu2007invitation}.

A \textit{stratification} of an $m$-dimensional variety $X\subseteq \CC^n$ is a filtration
\begin{align}
    \emptyset=F_{-1}\subseteq F_0\subseteq F_1\subseteq \cdots \subseteq F_m = X,
\end{align}
where $F_k$ are varieties, each $X_k:=F_k\setminus F_{k-1}$ is smooth and either empty or all its irreducible components are $k$-dimensional, and $\overline{X_k}\setminus X_k\subseteq F_{k-1}$. The stratification satisfies \textit{Whitney condition} (a) if, for every $0 \le
j<k\le m$, the pair $(X_j , X_k)$ satisfies: given a sequence $x_n\in X_k$ that converges in Euclidean topology to some $x\in X_j$ such that the tangent spaces $T_{x_n}\:X_k$ converge to a subspace $T$, we have that $T_x\:X_j\subseteq T$. A \textit{Whitney stratification} satisfies Whitney condition (a) and a stronger condition (b) that we do not need for this paper. By a theorem of Whitney, any algebraic variety has a Whitney stratification~\cite{whitney1992tangents,parusinski1995chern}.

We finally assert that, for a dominant morphism, the generic points of all irreducible components of a generic fiber are smooth points of the domain.

\begin{lemma}\label{le: gen fib} Let $f:X\to Y$ be a dominant polynomial map between irreducible affine varieties. There is no proper subvariety $Z\subseteq X$ that contains an irreducible component of the fiber $f^{-1}(y)$ for generic $y\in Y$. 
\end{lemma}

\begin{proof} Assume that $Z\subseteq X$ is a variety that contains a component of $f^{-1}(y)$ for generic $y\in Y$. Each component of such fibers is of dimension $\dim X-\dim Y$~\cite[Chapter 1, Section 8, Corollary 1]{mumford1999red}. Since $f^{-1}(y_1)\cap f^{-1}(y_2)=\emptyset$ for $y_1\neq y_2$, we conclude that $Z$ must be of dimension $(\dim X-\dim Y)+\dim Y=\dim X$. Since $X$ is irreducible, $Z=X$. 
\end{proof}

%%%%%%%%%%%%%%%%%%%%%%%%%%%%%%%%%%%%%%%%%%%%%%%%%%%%%%%%%%%%%%%%%%%%%%%%%%%%%%%%%%%%%%%%%%%%%%%%%%%%%%%%%%%%%%%%%%%%%%%%%%%%%%%%%%%%%%%%%%%%%%%%%%%%%%%%%%%%%%%%%%%%%%%%%%%%%%%%%%%%%%%%%%%%%%%%%%%%%

\section{\texorpdfstring{Computing $\mathrm{EDdeg}(X)$, $\Sigma_X$, and $\Sigma_X^\infty$}{}}\label{s: Comp}

In this section, we first recall how to compute the ED degree and the ED discriminant of a variety. We then present pseudo code in \Cref{alg: sigma} that computes the infinite Euclidean distance discriminant. This pseudo code is based on ideas that were kindly conveyed to us by Professor Mike Stillman at Cornell University. We prove its correctness in \Cref{thm: correctness}. 

Computing Euclidean distance degrees exactly is, in general, a difficult problem. Numerically, one may fix a data point \new{$u$ and} solve the critical equations; if $u$ is generic enough, we get exactly ED degree many points. However, because we do not a priori know if $u$ is generic enough, such numerical computations only provide conjectural values. For certified determination of ED degrees, there are analytical approaches based on Morse theory~\cite{EDDegree_point} and a computational approach using the ED polynomial. The following result is known (although a slightly incomplete proof appears in~\cite{horobect2019offset}); we include a proof for completeness:

\begin{theorem}\label{thm: deg EDpoly} Let $X$ be an irreducible variety in $\CC^n$ and assume $\mathrm{EDdeg}(X)>0$. The ED polynomial (defined over $\CC$) is of degree $2\cdot \mathrm{EDdeg}(X)$ in $\epsilon$. 
\end{theorem}

The assumption that $\mathrm{EDdeg}(X)>0$ is made to ensure that the ED polynomial is \new{well-defined. A} sufficient condition for $\mathrm{EDdeg}(X)>0$ is given in~\cite[Theorem 4.1]{draisma2016euclidean}: namely, that $T_x\:X \cap (T_x\:X)^\bot=\{0\}$ at some smooth $x\in X$. This holds, for instance, if $X$ is cut out by polynomials with real coefficients and its real dimension coincides with its complex dimension.

\begin{remark}\label{re: cond} See \cite[Section 7 \& 8]{catanese2000focal} for more details on varieties with ED degree 0 in the projective setting. 
\end{remark}

\begin{proof} We prove that the ED degree is $\deg_{R} \: F_u(R)$, where $F_u(R)$ is \new{the ED polynomial as defined in} \Cref{s: prelim}. Since the bisector locus is at most a hypersurface (\Cref{prop: bisec}), \new{all squared distances $(u-x)\cdot (u-x)=R$ are distinct among $x\in \mathcal E_X(u)$, for generic $u$}. If $(x,u)$ lies in the ED correspondence, then $F_u(d_u(x))=0$ and therefore $\deg_{R}\:F_u(R)\ge \mathrm{EDdeg}(X)$.

By Chevalley's theorem~\cite[Theorem 3.16]{JoeHarris}, $\mathcal O_{X}$ is the Euclidean closure of the image of $(x,u)\mapsto (u,d_u(x))$ from $\mathcal E_X$. In particular, for generic $u$, there are no additional solutions in $R$ to $F_u(R)=0$ other than those for which there is an $x$ with $R=d_u(x)$. \new{We are left to exclude the possibility that $F_u(R)$ has double roots in $R$ for generic $u$. In this direction, consider the map $\mathrm{der}_R: \mathcal O_X\to \CC$, sending $(u,R)\mapsto \frac{\partial}{\partial R} F_u(R)$. Assume first that $\mathrm{der}_R$ is dominant and consider the preimage $\mathrm{der}_R^{-1}(0)$. This preimage is a proper subvariety of $\mathcal O_{X}$, at most $(n-1)$-dimensional. Its projection onto the first factor $\CC^n$ is therefore at most $(n-1)$-dimensional, meaning that a generic data point $u$ has $\frac{\partial}{\partial R} F_u(R)\neq 0$ for every $R$ with $(u,R)\in \mathcal O_X$. This implies that there generically are no double roots in $R$. Assume second that $\mathrm{der}_R$ is not dominant. We only need to consider the case where $\mathrm{der}_R(\mathcal O_X)=\{0\}$. By locally parameterizing $\mathcal O_X$ via $(u,R(u))$, one can check that the normal of $\mathcal O_X$ is spanned by $(\frac{\partial}{\partial u_1}R(u),\ldots,\frac{\partial}{\partial u_n}R(u), -1)$. This however contradicts the assumption that $\mathrm{der}_R(\mathcal O_X)=\{0\}$.} 
\end{proof}

In the \new{source folder of the} \texttt{arXiv} version of this paper, we provide \texttt{Macaulay2} \new{code} that computes ED degrees (using \Cref{thm: deg EDpoly}), ED discriminants, bisector loci, and offset discriminants for varieties that are cut out by polynomials with rational coefficients. The correctness of the code follows from implicitization and the Elimination and Closure Theorems of~\cite[Chapter 3]{cox1997ideals}. To compute ED discriminants, we follow the discussion on the Implicit Function Theorem from \Cref{ss: local par}. Let $I$ be the prime ideal of $\mathcal E_X$, and let $I_u$ be the specialization ideal $\{f(x,u)\in \CC[x]: f\in I\}$. Then
\begin{align}
    \Sigma_X=\overline{\{u\in \CC^n: \mathrm{rank} \:\mathrm{Jac}_x(I_u)<n\textnormal{ for some }x\textnormal{ with }(x,u)\in \mathcal E_X\}}.
\end{align}
The bisector locus of $X$ is computed by first considering the ideal $I\subseteq \QQ[x,y,u]$ that is the sum of the prime ideal of $\mathcal E_X$ in $(x,u)$ \new{ and $(y,u)$ coordinates, as well as the ideal generated by the equidistance condition}. After saturating away the locus $x=y$, we eliminate the $x$ and $y$-variables to obtain the bisector locus.

For the computation of infinite ED discriminants, we make use of Gröbner bases~\cite{cox1997ideals}. We recall the basic definitions. Given a monomial ordering and an ideal $I\subseteq k[x]$ over a field $k$, a \textit{Gröbner basis} $G\subseteq k[x]$ is a finite set generating $I$ such that 
\begin{align}
   \langle \mathrm{LT}(g):g\in G\rangle  =  \langle \mathrm{LT}(f):f\in I\rangle,
\end{align}
where $\mathrm{LT}(f)$ is the leading term with respect to the monomial ordering. A \textit{reduced} Gröbner basis $G$ is one in which the leading coefficient of each polynomial is $1$, and for any $g\in G$, no monomial of $g$ lies in $ \langle \mathrm{LT}(h):h\in G\setminus \{g\}\rangle $. With this notation in place, we recall two results on Gröbner bases that we need for our algorithm. 

\begin{proposition}\label{prop: Finiteness} Let $G$ be a Gröbner basis for an ideal $I\subseteq \QQ[x]$. Then over $\CC$, the zero locus $V(I)$ is finite if and only if for each variable $x_i$, there is a $g\in G$ whose leading term is a power of $x_i$.  
\end{proposition}

\begin{proof} By~\cite[Proposition 3.7.1]{kreuzer2000computational}, $V(I)$ is finite if and only if for each variable $x_i$, a power of $x_i$ lies in the ideal of leading terms of $I$. By the equivalence of (i) and (ii) of~\cite[\S 5.3 Theorem 6]{cox1997ideals}, the statement follows. (Note that~\cite[\S 5.3 Theorem 6]{cox1997ideals} makes a subtly different finiteness claim.)
\end{proof}

We consider product monomial orders $\succ$ on $\CC[x,u]$ with the $x$-block preceding the $u$-block, i.e., exponents of the $x$-variables are compared first, and only in case of a tie are the $u$-variables compared. For our algorithm to be efficient, it is important to make a good choice of $\succ$. We choose the GRevLex block order.

\begin{proposition}\label{prop: specialization} Let $G$ be a Gröbner basis for an ideal $I \subseteq \QQ[x, u]$ with respect to a product monomial order with the $x$-block preceding the $u$-block. For each $g\in G\setminus \QQ[u]$, write
\begin{align}
    g = c_g(u)x^{\alpha_g} + (\textnormal{ terms }\prec x^{\alpha_g}).
\end{align}
If $b\in V(I\cap \QQ[u])$ and $c_g(b)\neq 0$ for all $g\in G\setminus \QQ[u]$, then 
\begin{align}
    \{g(x,b): g\in G\setminus \QQ[u]\}\subseteq \QQ[x]
\end{align}
is a Gröbner basis for the specialization ideal $I_b$.
\end{proposition}

\begin{proof} See~\cite[\S 4.7 Theorem 2]{cox1997ideals}. 
\end{proof}

With these results in mind, we present \new{\Cref{alg: sigma} for computing infinite ED discriminants}. In the \texttt{arXiv} version of this paper, we attach a \texttt{Macaulay2} implementation of this algorithm. Nothing in this algorithm is, in principle, specific to the nearest point problem; it can be adapted to compute the locus of $u\in \CC^m$ with positive-dimensional fibers under the projection onto the second factor for any algebraic variety $\mathcal E\subseteq \CC^n \times \CC^m$.

 \begin{algorithm}
\SetKwInOut{Input}{Input}\SetKwInOut{Output}{Output}
\SetKwInOut{Return}{Return}
\caption{Computation of $\Sigma_X^\infty$}\label{alg: sigma}
\Indm 
\Input{A proper subvariety $X\subseteq \CC^n$ cut out by polynomials with rational coefficients, and the GRevLex block monomial order on $\QQ[x,u]$, where the $x$-block precedes the $u$-block.}
\Output{The infinite Euclidean distance discriminant $\Sigma_X^\infty$.} 
  \Indp
    $I_\mathrm{EDC} \gets \textnormal{ Ideal of ED correspondence of }X\textnormal{ in }\QQ[x,u]$\;
      $\textnormal{\new{current\_ideals}}\gets \textnormal{ Set of one element } \langle 0\rangle \subseteq \QQ[x,u]$\;
       $\textnormal{inf\_ED\_ideals},\;\textnormal{new\_ideals}\gets \textnormal{ Empty set}$\;
     \While{$\textnormal{\new{current\_ideals}}\neq \emptyset $}{
       \For{$J$ in $\textnormal{\new{current\_ideals}}$}{
       \For{$K$ among the prime components of $(I_\mathrm{EDC} +J)\cap \QQ[u]$}{
         $G \gets \textnormal{Reduced Gröbner basis for }I_\mathrm{EDC} +K$\;
     $c_g(u)x^{\alpha_g}\gets \textnormal{ Leading term of }g\in G$\;
     \uIf{for each $i$, some $x^{\alpha_g}$ is a power of $x_i$ with $c_g(u)\in \QQ$}{Discard $K$\;}
    \uElseIf{for some $i$, no $x^{\alpha_g}$ is a power of $x_i$}{Add $K$ to $\textnormal{inf\_ED\_ideals}$\;}\Else{
    Add $K+\langle c_g(u)\rangle$ to $\textnormal{new\_ideals}$ for each $g\in G$ with $c_g(u)\not\in K$\;
    }
       }
    }
   
   $\textnormal{\new{current\_ideals}}\gets \textnormal{new\_ideals}$\;
     $\textnormal{new\_ideals}\gets \textnormal{ Empty set}$\;  }
       $I_\infty\gets \textnormal{ The intersection of all ideals in \textnormal{inf\_ED\_ideals}}$\;
\Indm 
\Return{$\textnormal{ The zero locus of }I_\infty$.}
 \end{algorithm}

\begin{theorem}[Code for infinite ED discriminants]\label{thm: correctness} \Cref{alg: sigma} terminates and correctly returns the infinite Euclidean distance discriminant.
\end{theorem}

\begin{proof} Since $X$ is cut out by polynomials of rational coefficients, this also holds for the ED correspondence\new{. This follows from implicitization theory \cite[Section 3.3]{cox1997ideals}.} 

For an ideal $J\subseteq \QQ[u]$, we proceed as follows. Let $K$ be a component of $(I_{\mathrm{EDC}}+J)\cap \QQ[u]$. Let $G$ be a reduced Gröbner basis of $I_{\mathrm{EDC}}+K$ for some product order. We consider the three possibilities:
\begin{enumerate}[label=(\Alph*)]
    \item For each $i$, there is a $g\in G$ whose leading $x$-monomial is a power of $x_i$ and whose coefficient in $u$ is constant.\label{enum: top consts}
    \item For some $x_i$, no leading $x$-monomial of $G$ is a power of $x_i$.\label{enum: no top}
    \item None of the above.\label{enum: else}
\end{enumerate}

If \Cref{enum: top consts} holds, then take a $u\in V(K)$ and consider the specialization ideal $I_{K,u}:=\{f(x,u)\in \QQ[x]: f\in I_{\mathrm{EDC}}+K\}$. It is spanned by $g(x,u)\in \QQ[x]$ for $g\in G$, and these $g(x,u)$ can be assumed to be contained in a Gröbner basis $G'$ for $I_{K,u}$. By \Cref{prop: Finiteness}, $V(I_K(u))$ cannot be infinite, meaning no point of the zero locus $V(K)$ lies in $\Sigma_X^\infty$.

If \Cref{enum: no top} holds, then no leading coefficient $c_g(u)$ for $g\in G\setminus \QQ[u]$ is in $K$, since $G$ is reduced. All such top coefficients $c_g(u)$ are nonzero for generic $u\in V(K)$, since $K$ is prime. We recover a Gröbner basis when specializing to such $u$, by \Cref{prop: specialization}. By assumption, there is some $x_i$ such that no leading $x$-monomial of $G$ is a power of $x_i$, and hence \Cref{prop: Finiteness} says that $V(I_{K,u})$ is at least one-dimensional. It follows that $V(K)$ is contained in $\Sigma_X^\infty$. 

If \Cref{enum: else} holds, then for generic $u\in V(K)$, the zero locus $V(I_{K,u})$ is finite. The only $u\in V(K)$ for which $V(I_{K,u})$ can be one-dimensional are those for which a top $x$-coefficient $c_g(u)$ for some $g\in G\setminus \QQ[u]$ vanishes, by \Cref{prop: specialization,prop: Finiteness}. \new{Therefore, we need to redo the case analysis for each ideal $K'=K+\langle c_g(u)\rangle$.} 

To see that the algorithm terminates, we consider $K'$ from \Cref{enum: else}. By the fact that $G$ is reduced, $K'$ is strictly bigger than $K$. Further, prime component $K''$ of $(I_{\mathrm{EDC}}+K')\cap \QQ[u]$ contains $K'$. This means that the dimension of $V(K'')$ is strictly smaller than that of $V(K)$.

\new{Finally,} we observe that in the first iteration of the while-loop, \new{current\_ideals} contains only $\langle 0\rangle$. The only prime component $K$ of $I_{\mathrm{EDC}}\cap \QQ[u]$ is the defining ideal of the closure of the image of $(x,u)\mapsto u$ restricted to $\mathcal E_X$. Since the infinite ED discriminant is contained in $V(K)$, \new{the algorithm is guaranteed to output} 
$\Sigma_X^\infty$. 
\end{proof}

%%%%%%%%%%%%%%%%%%%%%%%%%%%%%%%%%%%%%%%%%%%%%%%%%%%%%%%%%%%%%%%%%%%%%%%%%%%%%%%%%%%%%%%%%%%%%%%%%%%%%%%%%%%%%%%%%%%%%%%%%%%%

\section{\texorpdfstring{Examples of varieties with nonempty $\Sigma_X^\infty$}{}}\label{s: examples} 
\begin{figure}[ht]
\centering

\begin{minipage}{0.2\linewidth}
  \centering
  \includegraphics[width=\linewidth, trim=20 20 20 20, clip]{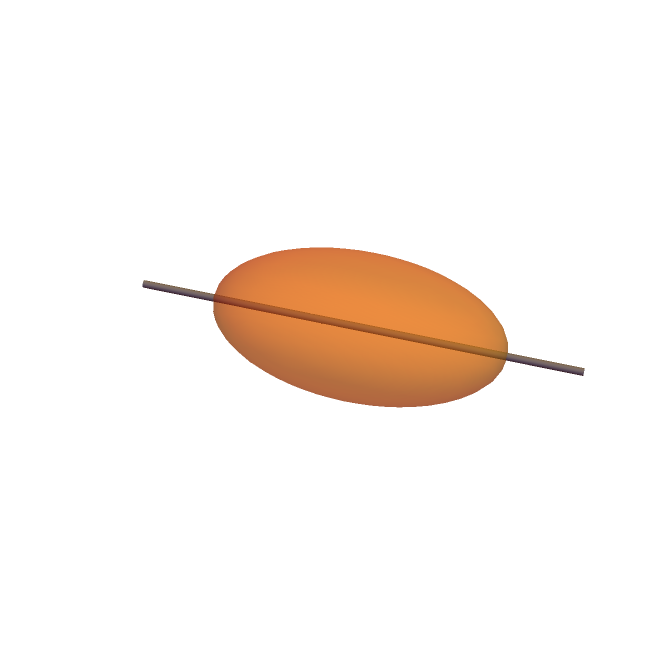} 
\end{minipage}\hspace{3em}
\begin{minipage}{0.2\linewidth}
  \centering
  \includegraphics[width=\linewidth, trim=20 20 20 20, clip]{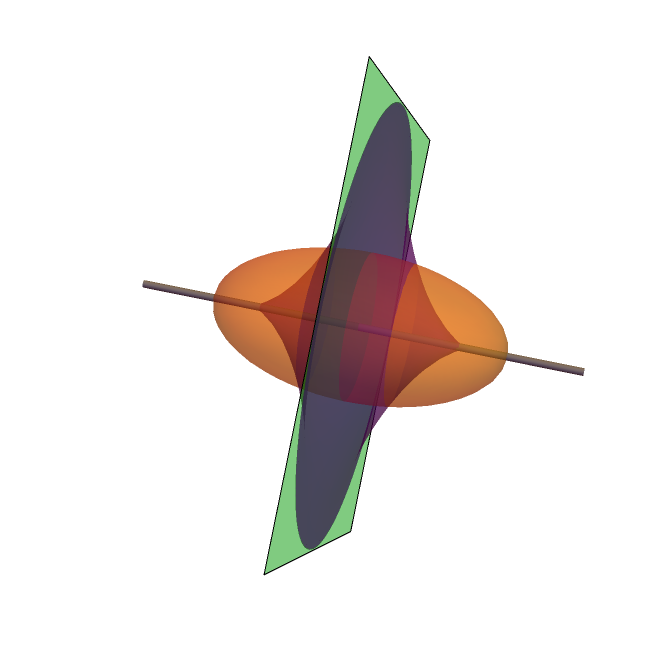} 
\end{minipage}\hspace{3em}
\begin{minipage}{0.2\linewidth}
  \centering
  \includegraphics[width=\linewidth, trim=20 20 20 20, clip]{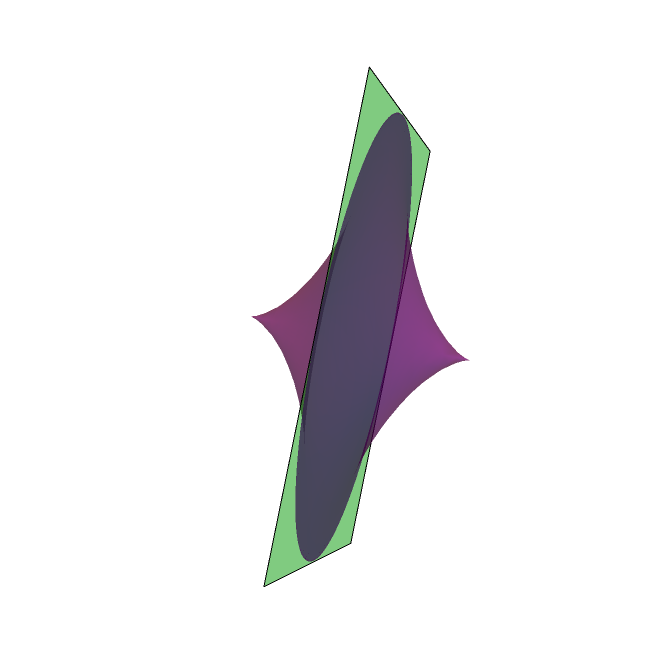} 
\end{minipage}

\vspace{-0.15cm}

\begin{minipage}{0.2\linewidth}
  \centering
  \includegraphics[width=\linewidth, trim=20 20 20 20, clip]{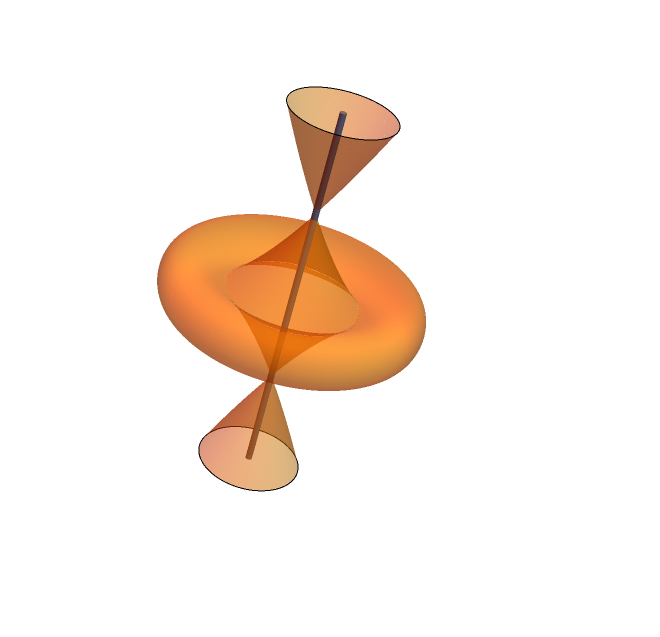} 
\end{minipage}\hspace{3em}
\begin{minipage}{0.2\linewidth}
  \centering
  \includegraphics[width=\linewidth, trim=20 20 20 20, clip]{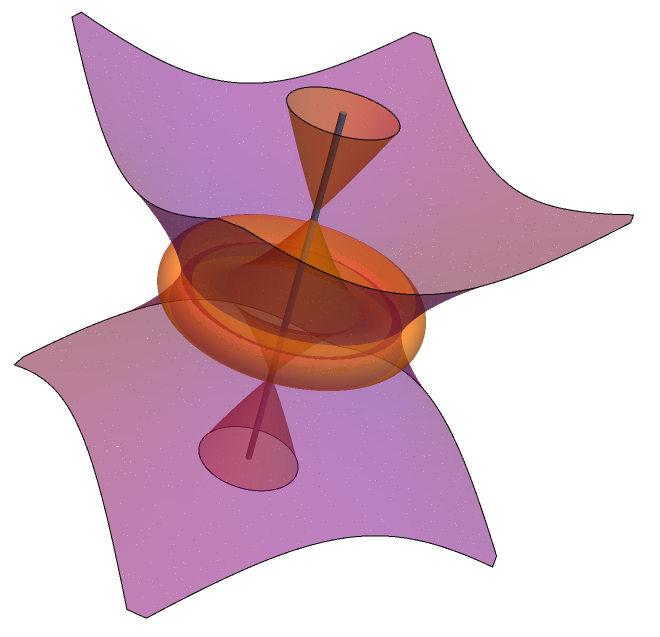} 
\end{minipage}\hspace{3em}
\begin{minipage}{0.2\linewidth}
  \centering
  \includegraphics[width=\linewidth, trim=20 20 20 20, clip]{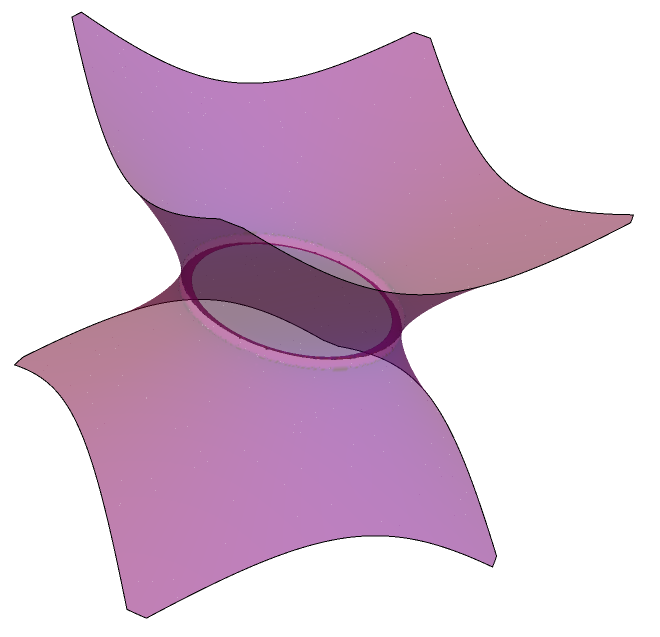} 
\end{minipage}

\vspace{-0.1cm}

\begin{minipage}{0.2\linewidth}
  \centering
  \includegraphics[width=\linewidth, trim=20 20 20 20, clip]{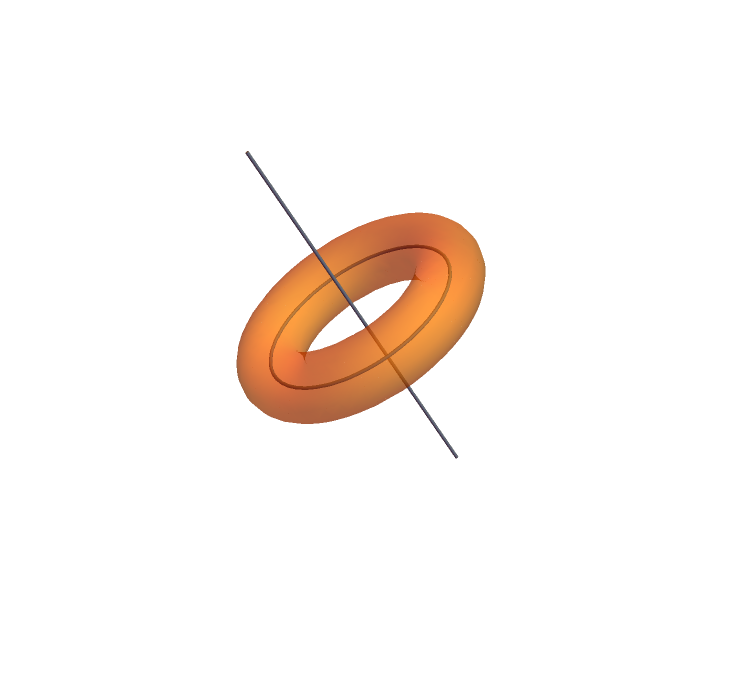} 
\end{minipage}\hspace{3em}
\begin{minipage}{0.2\linewidth}
  \centering
  \includegraphics[width=\linewidth, trim=20 20 20 20, clip]{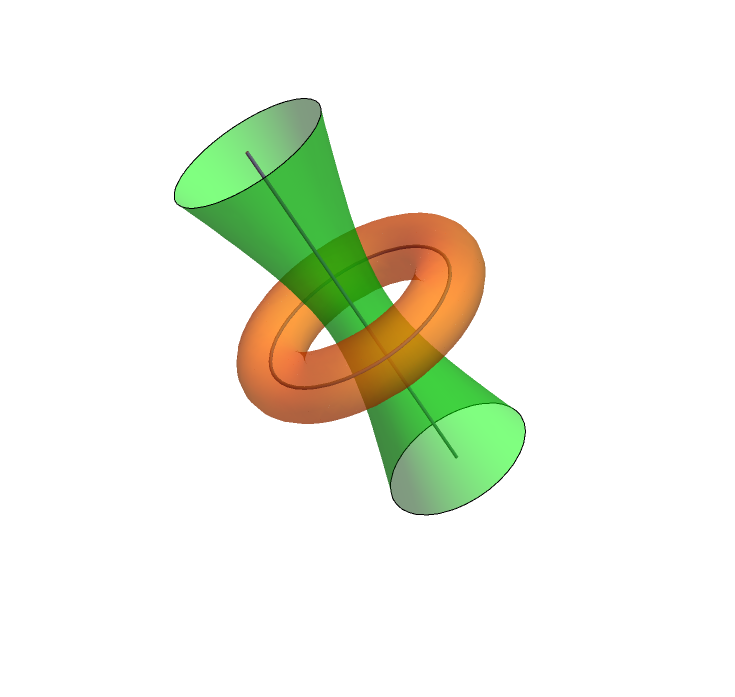} 
\end{minipage}\hspace{3em}
\begin{minipage}{0.2\linewidth}
  \centering
  \includegraphics[width=\linewidth, trim=20 20 20 20, clip]{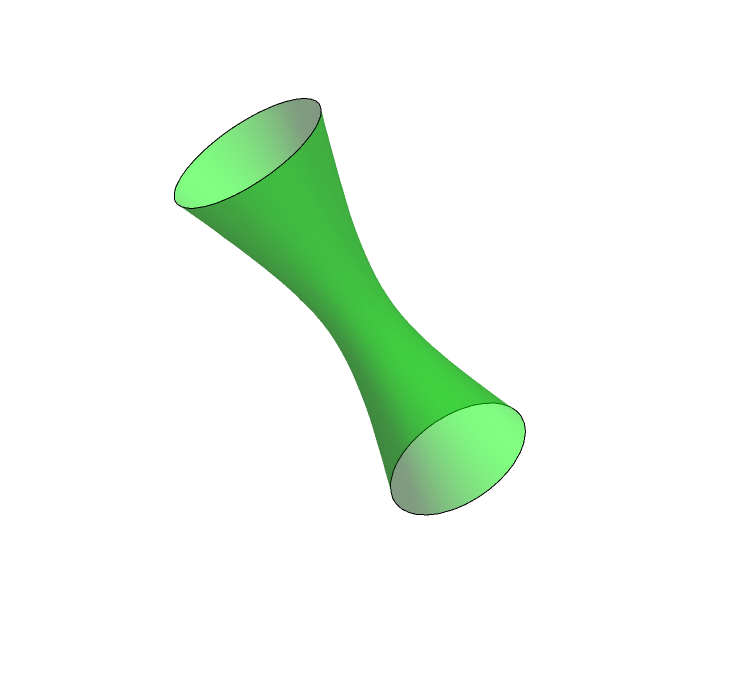} 
\end{minipage}

\vspace{-0.45cm}

\begin{minipage}{0.2\linewidth}
  \centering
  \includegraphics[width=\linewidth, trim=20 20 20 20, clip]{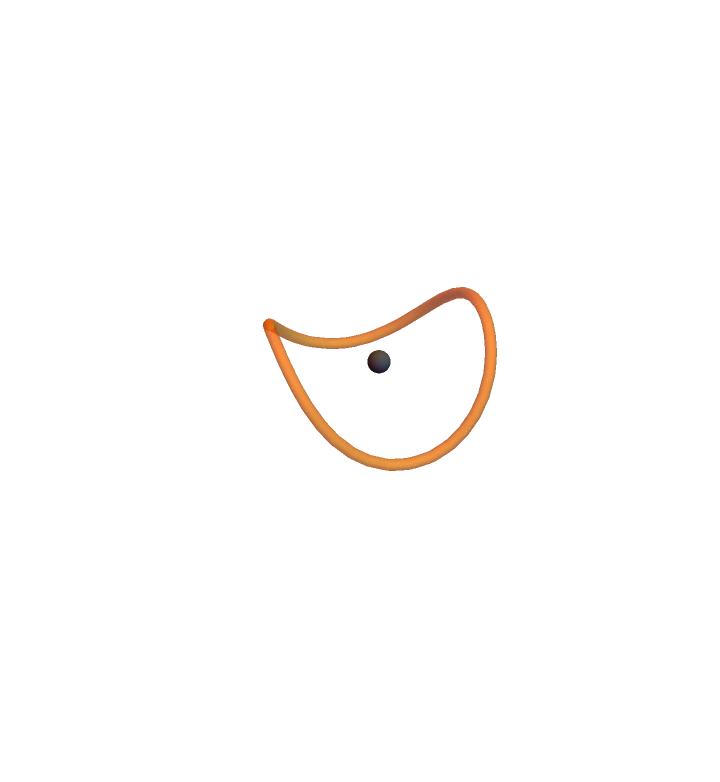} 
\end{minipage}\hspace{3em}
\begin{minipage}{0.2\linewidth}
  \centering
  \includegraphics[width=\linewidth, trim=20 20 20 20, clip]{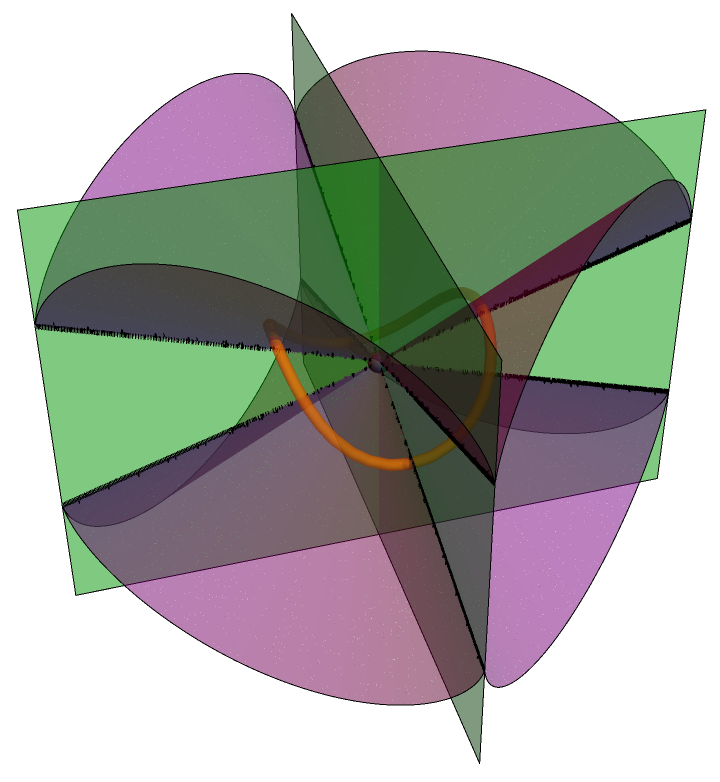} 
\end{minipage}\hspace{3em}
\begin{minipage}{0.2\linewidth}
  \centering
  \includegraphics[width=\linewidth, trim=20 20 20 20, clip]{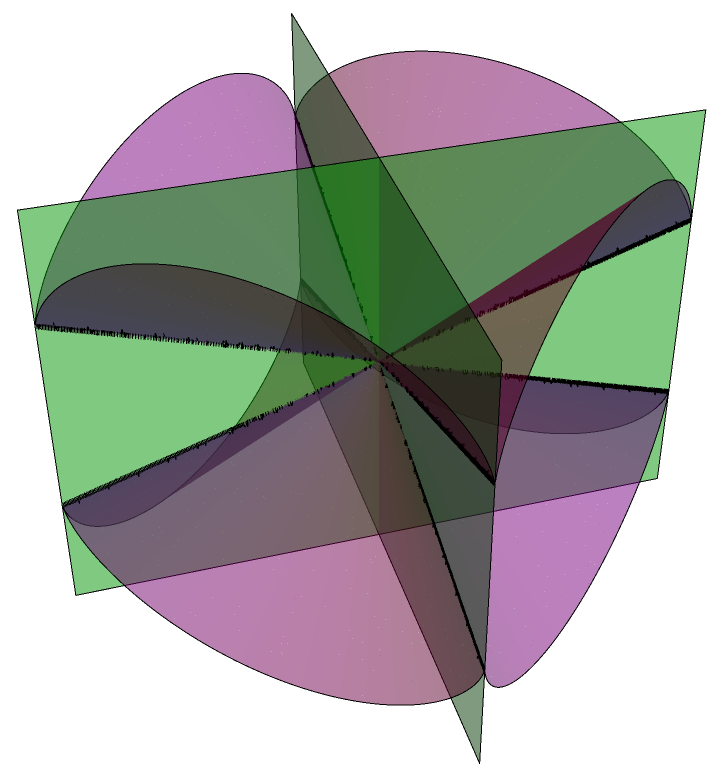} 
\end{minipage}

\caption{Varieties $X$ (orange) in $\CC^3$, and components of the \new{infinite ED discriminants (gray), ED discriminants (purple), and bisector loci (green). First column: $X$ and $\Sigma_X^\infty$. Second column: the first and third columns overlapped. Third column: real surface components of $\Sigma_X$ and $B_X$. The varieties displayed along the rows are as follows. First row: American football (\Cref{ex: foot}). Second row: volcano (\Cref{ex: volcano}). Third row: donut (\Cref{ex: donut}). Fourth row: potato chip (\Cref{ex: potato}).}}
\label{fig: examples}
\end{figure}

In this section, we present examples of varieties $X$ with nonempty infinite Euclidean distance discriminant. \new{Using \Cref{s: Comp}}, we \new{compute ED degrees $\mathrm{EDdeg}(X)$, ED discriminants $\Sigma_X$, bisector loci $B_X$, offset discriminants $\Delta_X$, and infinite ED distance discriminants $\Sigma_X^\infty$.} \new{The} code is run on \new{a} 13th Gen Intel(R) Core(TM) i9-13900H CPU
running at 2.60 GH. The first four examples \new{concern surfaces} of revolution, a tube, and a spherical curve \new{(a curve inside a sphere)} in $\CC^3$, illustrated in \Cref{fig: examples}. For their Euclidean distance correspondences in $\CC^3\times \CC^3$, we use coordinates $((x,y,z),(u,v,w))$. A general result on infinite ED discriminant for these types of varieties is given in \Cref{thm: rev and tubes}. 

The isotropic quadric \new{(\Cref{ex: isotropic})} plays a special role in nearest point problems. As long as $n\ge 3$, we show that its infinite ED discriminant coincides with the quadric itself. \new{We further examine infinite ED discriminants for multiview varieties in the case of two views, noting that the original motivation for the studying ED degrees arose from triangulation in computer vision \cite{hartley1997triangulation,stewenius2005hard}; see \cite[Example 3.3]{draisma2016euclidean}.} In \Cref{ex: epip}, we show that for a generic fundamental matrix $F^{12}$, $\Sigma_X^\infty$ is \new{empty and} we describe those $F^{12}$ for which it is nonempty. In \Cref{prop: linear}, we determine precisely which affine linear varieties have nonempty $\Sigma_X^\infty$. Another application for the nearest point problem is low-rank matrix approximation, where the infinite ED discriminant is of codimension 2 (\Cref{prop: rank-r}).

\subsection{First examples} 

\new{Surfaces of revolution generally have the property that the axes they revolve around are contained in their infinite ED discriminants. Such a surface is a union of circles, and these circles are sets of critical points for data on the axis. A detailed argument is provided in \Cref{ss: five}. Our first three examples, an American football, a volcano, and a donut, are surfaces of revolution with this property.}

\begin{example}[American football]\label{ex: foot}
Let $X\subseteq \CC^3$ be the symmetric ellipsoid defined by the equation 
\begin{align}
    x^2+4y^2+4z^2=4.
\end{align}
The ED degree is $4$. \new{The infinite ED discriminant is the $u$-axis. The ED discriminant is the union of the $u$-axis and a Lam\'e curve (the ED discriminant of the ellipse $x^2+4z^2=0$) revolved around the $u$-axis. The bisector locus consists of 2 components: the $u$-axis and the $vw$-plane. Our computations show that $\Delta_X$ is not the union $\Sigma_X\cup B_X$; the offset discriminant contains the additional union of hyperplanes $v^2+w^2=0$.}
\hfill$\new{\rotsymbol\,}$
\end{example}

\begin{figure}[ht]
\centering

\begin{minipage}{0.5\linewidth}
  \centering
  \includegraphics[width=\linewidth]{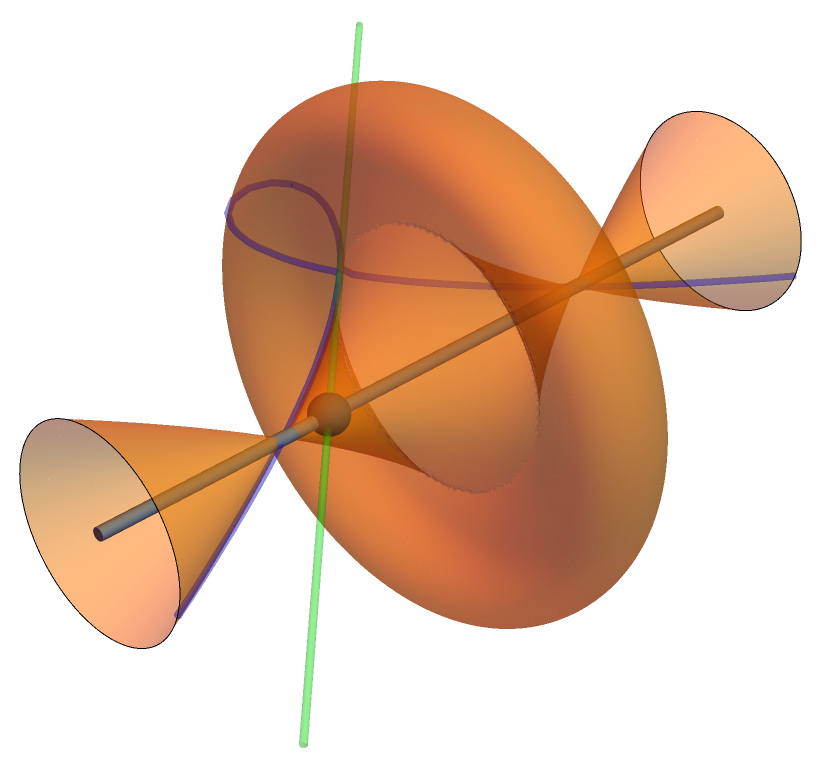} 
\end{minipage}

\caption{The volcano surface (orange), the nodal cubic (blue), the $x$-axis (gray), the data point $(1,0,0)$ (gray), and a (in the limit) normal line \new{through the singularity} (green) to the nodal cubic.}
\label{fig: volcano-cont}
\end{figure}

\begin{example}[Volcano]\label{ex: volcano} Let $X$ be the nodal cubic plane curve
\begin{align}\label{eq: plane cubic}
    x^2=(y+1)^2(y+2),\quad z=0
\end{align}
revolved around the $x$-axis. It is given by:
\begin{align}\begin{aligned}
    (y^2+z^2)^3-6(y^{2}+z^2)^2+(8x^{2}+9)(y^{2}+z^{2})-(x^{2}-2)^2=0.
\end{aligned}
\end{align}
This surface resembles two cone shaped volcano with smoke coming out at the top. The ED degree is $10$. \new{The infinite ED discriminant consist of five components: the $u$-axis and the union of four lines $u^2=2$, $v^2+w^2=0$. The ED discriminant additionally contains a degree-12 surface defined by a polynomial of 62 terms, the circle $u=0$, $v^2+w^2=1$, and the circles $u^2+18=0$, $v^2+w^2+5=0$. The computation of the bisector locus failed to terminate with our code.

Among the ED degree, the bisector locus, and the discriminants, $\Sigma_X^\infty$ was the fastest to compute, taking 6 seconds. In comparison, $\mathrm{EDdeg}(X)$ took $3$ minutes, $\Sigma_X$ took $4$ minutes, $B_X$ failed to terminate, and $\Delta_X$ required $15$ hours (one of its components is a degree-38 surface defined by a polynomial of 1166 term). 
}

\new{This volcano variety exhibits that $C_X(u)$ is not always dense in $\mathcal E_X(u)$. In this direction, note that one of the two normal lines through the nodal singularity $x=(0,-1, 0)$ meets the $x$-axis in $u=(1,0,0)$. By taking a converging sequence of data points $u_n\to u$ along the $x$-axis with associated critical points $x_n\to x$, we see that $\mathcal E_X(u)$ contains the circle $x=0, y^2+z^2= 1$. Our \texttt{Macaulay2} computations show that this containment is, in fact, an equality. However, since this circle lies in the singular locus of $X$, we have $C_X(u)=\emptyset$. We illustrate this phenomenon in \Cref{fig: volcano-cont}.} %https://www.desmos.com/3d/tl4dviq6c6
\hfill$\new{\rotsymbol\,}$\end{example}

\new{The next example, the donut, is also a surface of revolution. Additionally, it is a \textit{tube}, meaning an $\epsilon$-offset hypersurface associated to a curve. By construction, a tube is a union of spheres centered at points on the curve. The points on those spheres are all critical to those points on the curve. This follows from \Cref{le: offset duality} and is proven in \Cref{ss: five}.}

\begin{example}[Donut]\label{ex: donut} Let $X$ be the torus around the central circle $x^2+y^2=9,z=0$ of radius $1$. Over the real numbers, it satisfies $(3-\sqrt{x^2+y^2})^2+z^2=1$, which we expand to the algebraic equation
\begin{align}
    (x^2+y^2+z^2+3^2-1^2)^2-4\cdot 3^2\cdot(x^2+y^2)=0.
\end{align}
The ED degree is $4$. \new{The infinite ED discriminant coincides with the ED discriminant and consists of 3 components: the circle $w=0$, $u^2+v^2=9$, the $w$-axis, and the union of four lines $w^2+8=0$, $u^2+v^2=0$. The bisector locus is the union of the first two $\Sigma_X^\infty$ components and the surface $w^2-8u^2-8v^2+8=0$. For this variety, $\Sigma_X^\infty$ coincides with the intersection $\Sigma_X\cap B_X$.}
%https://www.desmos.com/3d/a2fvfywwoy
\hfill$\new{\rotsymbol\,}$
\end{example}

\new{The above examples raise the question: are there} surfaces in $\CC^3$ with a one-dimensional infinite discriminant other than tubes and surfaces of revolution? The answer is yes. In \Cref{s: skew}, we define \textit{skew-tubes} and show how to construct surfaces with this property. 

\new{Spherical curves $X$ exhibit nonempty infinite ED discriminants. Indeed, the affine normal space $x + N_x\: X$ contains the center of the sphere $c$ for every smooth $x\in X$. This fact is formally proven in \Cref{cor: center hyp}.}

\begin{example}[Potato chip]\label{ex: potato} Let $X$ be the spherical curve
\begin{align}
    x^2+y^2+z^2=1,\quad z=xy.
\end{align}
This degree-$4$ curve is reminiscent of the contour of a potato chip. The ED degree is $8$. \new{The infinite ED discriminant is the origin. The ED discriminant is the zero locus of
\begin{align}\begin{aligned}\label{eq: potato disc}
    &729u^{10}v^2-2484u^8v^4+3574u^6v^6-2484u^4v^8+729u^2v^{10}+729u^{10}w^2+1053u^8v^2w^2-\\
    &2166u^6v^4w^2-2166u^4v^6w^2+1053u^2v^8w^2+729v^{10}w^2+3591u^8w^4-636u^6v^2w^4-\\
    &3606u^4v^4w^4-636u^2v^6w^4+3591v^8w^4+5536u^6w^6-4512u^4v^2w^6-4512u^2v^4w^6+\\
    &5536v^6w^6+384u^4w^8-9984u^2v^2w^8+384v^4w^8-6144u^2w^{10}-6144v^2w^{10}-4096w^{12}.
\end{aligned}
\end{align}
The bisector locus consists of 3 components: the planes $u-v=0$ and $u+v=0$, and the surface $(u^2+v^2+2w^2)^2+4u^2v^2=0$. In particular, $\Sigma_X^\infty$} is neither a component of $\Sigma_X$, $B_X$, nor $\Sigma_X\cap B_X$.
%https://www.desmos.com/3d/ekiedlldqf
%\hfill$\rotsymbol\,$
%https://www.desmos.com/3d/nnypqdkdz2
%https://www.desmos.com/3d/hed7lg4fsj
%https://www.desmos.com/3d/dj5kk4gxku
\hfill$\new{\rotsymbol\,}$
\end{example}

 \begin{example}[Isotropic quadric]\label{ex: isotropic} The isotropic quadric $Q$ defined by $x_1^2+\cdots+x_n^2=0$ in $\mathbb{C}^n$, $n\ge 3$, enjoys several unique optimization properties. The ED degree is $0$, since for smooth $x\in Q$, the normal space $\{x+v: v\in N_x\: Q\}$ is completely contained in $Q$. \new{For this variety, $\Sigma_Q$, $B_Q$, $\Delta_Q$, and $\Sigma_Q^{\infty}$ coincide} and consist of one component: the isotropic quadric itself. Indeed, for data points $u\in Q\setminus \{0\}$, all scalings $\lambda u\in Q$ are critical points for $u$. This is because $\mathcal{E}_Q$ is generated by $x_1^2+\cdots +x_n^2=0$ and $x_iu_j-x_ju_i=0$, for $i\neq j$. To conclude: if $u\in \CC^n\setminus Q$, then $\mathcal E_Q(u)=\emptyset$; if $u=0$, then $\overline{C_X(u)}$ is the whole $Q$; and for any other $u\in Q$, the set $C_X(u)$ is one-dimensional. Note also that the squared distance between $x$ and $u$ such that $(x,u)\in \mathcal E_Q$ is 0. 
\hfill$\new{\rotsymbol\,}$
\end{example}

We shift our focus to a variety that appears in computer vision~\cite{Hartley2004}. A camera is a $3\times 4$ full rank matrix $C$. Given two such cameras $C_1, C_2$ of distinct \textit{centers}, that is, kernels, we study the joint image map
\begin{align}\begin{aligned}
    \Phi: \PP^3&\dashrightarrow \PP^2\times \PP^2,\\
    P&\mapsto (C_1P,C_2P).
\end{aligned}
\end{align}
The image of $\Phi$ is the set of all possible pairs of pictures that can be taken of a point in $\PP^3$ by $C_1,C_2$, and its closure is a \textit{multiview variety}. This variety is cut out by a bilinear form 
\begin{align}
    x^\top F^{12}y=0, 
\end{align}
known as the \textit{epipolar constraint}, where $(x,y)\in \PP^2\times \PP^2$ are coordinates for the image points. The \textit{fundamental matrix} $F^{12}$ is a $3\times 3$ rank-$2$ matrix. Given a noisy image pair $(\widetilde{x},\widetilde{y})$ in an affine patch $\RR^2\times \RR^2\subseteq \PP^2\times \PP^2$, \textit{triangulation} refers to recovering a world point in an affine patch $\RR^3\subseteq \PP^3$ that best matches the image pair $(\widetilde{x},\widetilde{y})$; this is a nearest-point problem. 

\begin{example}[Epipolar constraint]\label{ex: epip} Let $F$ be a real $3\times 3$ rank-$2$ matrix. Define $X$ to be the variety in $\CC^4$ of points $((x_1,x_2),(x_3,x_4))\in \CC^2\times \CC^2$ satisfying
\begin{align}
    \begin{bmatrix}
        x_1 & x_2 & 1
    \end{bmatrix}F\begin{bmatrix}
        x_3\\ x_4 \\1
    \end{bmatrix}=0.
\end{align}
Assume that the top left $2\times 2$ submatrix $F_{2\times 2}$ of $F$ is invertible. The epipolar constraint is a degree-2 expression, and by translation and orthogonal action, it equals
\begin{align}\label{eq: w}
    y_1^2-y_2^2+wy_3^2-wy_4^2=0,
\end{align}
for some $w\neq 0$; the details are worked out in~\cite{rydell2025framework}. Using the tools of \Cref{s: Comp}, we have verified that as long as $w\neq\pm 1$, the infinite ED discriminant of \eqref{eq: w} is empty. For $w=\pm 1$, $\Sigma_X^\infty$ is a union of four planes. In the case that $w=1$, it consists of four components: \new{the planes $u_2=u_4=0$ and $u_1=u_3=0$,} the surface
\begin{align}
    u_2u_3-u_1u_4 =u_2^2+u_4^2=u_1u_2+u_3u_4=u_1^2+u_3^2=0,
\end{align}
which is the set of rank-1 matrices 
\begin{align}\label{eq: Umat}
   U= \begin{bmatrix}
        u_1 & u_2\\ u_3 & u_4
    \end{bmatrix}
\end{align}
for which $U^\top U=0$, and \new{the last component is this surface} after permuting $u_1$ with $u_3$. As proven in~\cite{rydell2025framework}, the ED degree drops to $2$ precisely when $w=\pm 1$, suggesting that a nonempty infinite ED discriminant is correlated with lower ED degrees. This is explored further in \Cref{ss: opt}.
\hfill$\new{\rotsymbol\,}$
\end{example}

A classic instance of nearest point problems is closest rank-$r$ fitting. Consider the variety of rank-$r$ $n\times m$-matrices $M_{n\times m}^{\leq r}$. The celebrated Eckart-Young theorem and~\cite[Example $2.3$]{draisma2016euclidean} shows for a generic real matrix $U$ with singular value decomposition $SDT^\top$, the critical points are 
\begin{align}
    S\cdot \text{diag}(0,\ldots,\sigma_{i_1},0,\ldots,0,\sigma_{i_r},0,\ldots,0)\cdot T^\top,
\end{align}
where $\sigma_{i}$ is the $i$-th largest singular value of $U$. In particular, 
\begin{align}\label{eq: deg rank-r}
    \mathrm{EDdeg}(M_{n\times m}^{\leq r}) = \binom{\min\{n,m\}}{r}.
\end{align}

\begin{example}[Cone over smooth quadric]\label{ex: cone over quadric} The equation  
\begin{align}
    x_1x_4-x_2x_3=0
\end{align}
defines a smooth projective quadric surface in $\PP^3$. Let $X$ be its cone over $\CC^4$. It is the variety of rank-deficient $2\times 2$ matrices
\begin{align}
    \begin{bmatrix}
        x_1 & x_2 \\ x_3 & x_4
    \end{bmatrix}.
\end{align}
The ED degree is 2, which agrees with the formula \eqref{eq: deg rank-r}. \new{The infinite ED discriminant coincides with the ED discriminant and the bisector locus, and consists of 4 components: the plane $u_2-u_3=u_1+u_4=0$, the plane $u_2+u_3=u_1-u_4=0$, the surface of rank-1 matrices \eqref{eq: Umat} for which $UU^\top=0$, and the surface of rank-1 matrices \eqref{eq: Umat} for which $U^\top U=0$. The first and second components are, up to scaling, the two components of matrices $U$ satisfying $U^\top U=I$.}
\hfill$\new{\rotsymbol\,}$
\end{example}

%%%%%%%%%%%%%%%%%%%%%%%%%%%%%%%%%%%%%%%%%%%%%%%%%%%%%%%%%%%%%%%%%%%%%%%%%%%%%%%%%%%%%%%%%%%%%%%%%%%%%%%%%%%%

\subsection{Five families of varieties}\label{ss: five}

We move on to \new{state general results on} infinite ED discriminants for five families of varieties: affine linear spaces, spherical varieties, $\epsilon$-offset hypersurfaces, varieties of revolution, and varieties of rank-$r$ \new{matrices.} 

\begin{proposition}\label{prop: linear} Consider \new{a nonempty} affine linear space $X=\{x\in \CC^n: Ax=b\}$, where $A$ is a full rank \new{$m\times n$ matrix $A$ with $m\le n$}. The discriminants $\Sigma_X,\Sigma_X^\infty$ coincide. They are nonempty if and only if $\det AA^\top=0$, and in this case, they \new{equal the preimage of~~$\mathrm{Im}\; AA^\top +b$ with respect to $A:\CC^n\to \CC^m$}, which is of dimension $\new{n-m}+\mathrm{rank}(AA^\top)$.
\end{proposition}

The ED degree of the affine space $X$ is one if and only if \new{$\mathrm{rank}(AA^\top)=m$}, and zero otherwise~\cite[Proposition 1.2.7]{sodomaco2020distance}. Do there exist other varieties with ED degree one? The answer turns out to be yes~\cite[Remark 3.7]{maxim2020defect}.

\begin{proof} Assume first that $AA^\top$ is invertible. \new{For} any data point $u$, the unique critical point is $x=u-A^\top (AA^\top)^{-1}(Au-b)$. This means that $\Sigma_X=\Sigma_X^\infty=\emptyset$. Second, assume that $AA^\top$ is singular. Now, all $u$ that have critical points have infinitely many. These $u$ are those such that $Au-b\in \mathrm{Im}\: AA^T$. This set makes up both $\Sigma_X$, $\Sigma_X^\infty$, and is \new{the preimage from the statement}. 
\end{proof}

Our next theorem captures \Cref{ex: foot,ex: volcano,ex: donut} as \new{special cases}. In this direction, we generalize surfaces of revolution as follows. Let $V\subseteq W\subseteq \CC^n$ be strict inclusions of affine linear spaces. Let $Y$ be an irreducible variety in $W$. For $y\in Y$, denote by $\mathrm{pr}_V(y)$ the orthogonal projection onto $V$. Let $N_W(y)$ denote the span of \new{$y- \mathrm{pr}_{V}(y)$ and the normal space of $W$}. The \textit{variety of revolution} $\mathrm{Rev}(V,W,Y)$ is the union of spheres
\begin{align}\label{eq: rev def}
    \bigcup_{y\in Y} \new{\Big\{}\mathrm{pr}_{V}(y) + b:\; b\in N_W(y),\; b\cdot b = (y- \mathrm{pr}_{V}(y))\cdot(y- \mathrm{pr}_{V}(y)) \new{\Big\}}. 
\end{align}
For simplicity, assume $V$ is spanned by the first $k$ unit vectors $e_1,\ldots,e_k$, and $W$ by $V$ and $e_{k+1},\ldots,e_m$. Then an alternative description of $\mathrm{Rev}(V,W,Y)$ is the closure of 
\begin{align}
\left\{\, x \in \CC^n \;:\;
\begin{aligned}
    &\textnormal{there exists a }\lambda\in\CC \textnormal{ such that}\\ &(x_1,\ldots,x_k,\lambda x_{k+1},\ldots,\lambda x_m\new{,0,\ldots,0})\in Y\\
    &\textnormal{and }\lambda^2\sum_{i={k+1}}^m x_i^2 = \sum_{i={k+1}}^nx_i^2
\end{aligned}
\right\}.
\end{align}

\begin{theorem}\label{thm: rev and tubes} $ $

\begin{enumerate}
\item A variety $X$ contained in a hypersphere $S$ of center $c$ has that $c\in \Sigma_X^\infty$.\label{enum: hyper}
    \item Assume that $Y\subseteq \CC^n$ is irreducible of dimension at most $n-2$ \new{with $\mathrm{EDdeg}(Y)>0$. Then $Y\subseteq \Sigma_X^\infty$ for some irreducible component $X$ of the $\epsilon$-offset hypersurface $\mathcal O_{Y,\epsilon}$.} % Any $\epsilon$-offset hypersurface $X=\mathcal O_{Y,\epsilon}$ has that $Y\subseteq \Sigma_X^\infty$.
    \label{enum: off}
    \item Assume that $V\subseteq W\subseteq \CC^n$ is a strict inclusion of affine linear spaces. Assume $Y\subseteq W$ is an irreducible variety such that $\mathrm{pr}_V:Y\to V$ is dominant. If for generic $u\in V$, $C_Y(u)\setminus V$ is nonempty, then the variety of revolution $X=\mathrm{Rev}(V,W,Y)$ has that $V\subseteq \Sigma_X^\infty$.\label{enum: rev}
\end{enumerate}
\end{theorem}

%\new{For \Cref{enum: off}, we point out that $X$ may be reducible, and in that case we mean that at least one component of $X$ has $Y$ in its infinite ED discriminant. It is also an implicit assumption that $\mathrm{EDdeg}(Y)>0$, because the ED polynomial is not defined otherwise.} 
Before the proof, let $T\: V$ denote the tangent space of an affine linear space $V$, and $N\: W$ the normal space of an affine linear space $W$. 

\begin{proof} $ $

For \Cref{enum: hyper}, we postpone the proof to \new{\Cref{cor: center hyp}}.

\new{For \Cref{enum: off}, consider} the rational dominant map 
\begin{align}\begin{aligned}
    Y\times \CC^{n-1-\new{\dim Y}}&\dashrightarrow X,\\
    (y,\new{t})&\mapsto y+\new{\phi}(y,\new{t}),
\end{aligned}
\end{align}
where \new{for fixed $y$, $\phi(y, t)$ is a local parametrization of the (at least one-dimensional) hypersphere inside $N_y\: Y$ centered at $0$ of squared radius $\epsilon^2$. In the case that $\dim N_y\: Y=2$ and $\epsilon =0$, the hypersphere is reducible, and we let $\phi(y,t)$ parametrize one such component instead. We can use \Cref{le: offset duality} to deduce that $Y\subseteq \Sigma_X^\infty$. Indeed, let $y(s) \in Y$ be a local parametrization and define $u(s, t) = y(s) + \phi(y(s), t)\in \mathcal O_{Y,\epsilon}$. Then the lemma says that $u(s, t)$ in $ \mathcal O_{Y,\epsilon}$ is critical to the data point $y(s)$ for each $t$.}

For \Cref{enum: rev}, fix a generic $u\in V$, a $y\in C_Y(u)\setminus V$, and some generic $x$ in the corresponding sphere $B(y)$ from \eqref{eq: rev def}. By construction, $x=O(y-v_0)+v_0$ for some orthogonal matrix $O$ that is the identity on $T\: V$, and some $v_0\in V$. We also have
\begin{align}
  \{O(y-v_0)+v_0: y\in Y\}\subseteq X
\end{align}
and \new{therefore the tangent space $T_x\:X$} contains $T_{Oy}\:OY =OT_y\: Y$. Moreover, $T_x\:X$ contains $T_x\:B(y)$, which equals $ON\:W$, because the tangent space of $B(y)$ at $y$ is $N\: W$. These tangent spaces together span $T_x\: X$ for dimension reasons. \new{To prove that $x$ is critical to $u$, we} need to check that for any $z\in T_x\: X$, we have $z\cdot (u-x)=0$. If $z\in OT_y\: Y$, then $z=O\hat z$ for some $\hat z\in T_y\: Y$. \new{It follows that $z\cdot (u-x)$} equals $\hat z \cdot O^\top (u-v_0+(Oy-v_0))$. Since $O$ is the identity on $T\: V$, we have $O^\top (u-v_0)=u-v_0$. As a consequence, $z\cdot (u-x)=\hat z \cdot (u-y)$, \new{which} is $0$ by assumption. An analogous calculation shows that if $z\in ON\: W$, then $z\cdot (u-x)=0$. We have shown that $B(y)\subseteq C_X(u)$, \new{which suffices since the sphere} $B(y)$ is at least of dimension $1$.   
\end{proof}

In \Cref{ex: cone over quadric}, we saw that the variety of rank-deficient $2\times 2$ matrices has a codimension 2 infinite ED discriminant. This is also true in general:

\begin{proposition}\label{prop: rank-r} Let $r<\min\{n,m\}$. The infinite ED discriminant of $M_{n\times m}^{\leq r}$ contains, as a component, the (codimension 2) complex closure of the variety of real matrices $U$ with two equal singular values. 
\end{proposition}

\begin{proof} We may without restriction assume \new{that} $n\le m$. It suffices to take a real full rank matrix $U$ with two equal singular values \new{and show} that $U\in \Sigma_X^\infty$, where $X=M_{n\times m}^{\leq r}$. Let $\sigma_1\ge \ldots\ge \sigma_n>0$ denote the singular values of $U$, and let $D$ denote the diagonal matrix with $\sigma_i$ at entry $(i,i)$. We follow the proof of~\cite[Theorem 2.9]{ottaviani2015geometric}. Any critical matrix $X$ to $U$ is given by 
\begin{align}
    S\cdot (D_{i_1}+\cdots +\new{D_{i_r}})\cdot T^\top,
\end{align}
where $i_j$ are distinct indices among $\{1,\ldots,n\}$\new{, $D_{k}$ is $\sigma_{k}$ at entry $(k,k)$ and zero everywhere else}, and $S,T$ are such that $U=SDT^\top$ is a singular value decomposition of $U$. \new{Let $s_i,t_i$ denote the columns of $S,T,$ respectively.}

Assume (without loss of generality) that $\sigma = \sigma_1=\sigma_2$. Then for the singular value decomposition of $U$, there are infinitely many choices of $s_{1},s_{2}$ and $t_{1},t_{2}$. Indeed, 
\begin{align}
    \begin{bmatrix}
        s_1 & s_2 
    \end{bmatrix} R \begin{bmatrix}
        \sigma & 0\\ 0& \sigma 
    \end{bmatrix}
    R^\top  \begin{bmatrix}
        t_1 & t_2 
    \end{bmatrix}^\top = \sigma s_1\new{t}_1^\top + \sigma \new{s}_2t_2^\top   
\end{align}
for any orthogonal $2\times 2$ matrix $R$. \new{Let} $s_1(R),s_2(R)$ be the first and second columns of $\begin{bmatrix}
        s_1 & s_2 
    \end{bmatrix} R$, and \new{define $t_1(R),t_2(R)$ analogously}. Then $(s_i(R),t_i(R))$ are pairs of singular vectors. In particular,
    \begin{align}
        \new{\sigma} s_2(R) \new{t}_2(R)^\top +\new{\sigma_3 s_{3}t_{3}^\top +}\cdots + \sigma_{r+1} s_{r+1}t_{r+1}^\top
    \end{align}
    is a critical point \new{to $U$} for every $R$. Choosing $R$ as the identity \new{or} the reflection along the $x=y$ axis yields two different critical points. By continuity \new{of $s_i(R),t_i(R)$}, there are infinitely many different ones. 
\end{proof}

%%%%%%%%%%%%%%%%%%%%%%%%%%%%%%%%%%%%%%%%%%%%%%%%%%%%%%%%%%%%%%%%%%%%%%%%%%%%%%%%%%%%%%%%%%%%%%%%%%%%%%%%%%%%%%%%%%%%%%%%%%%%%%%%%%%%%%%%%%%%%%%%%%%%%%%%%%%%%%%%%%%%%%%%%%%%%%%%%%%%%%%%%%%%%%%%%%%%%

\section{Properties of infinite ED discriminants}\label{s: Properties}

In this section, we prove in full generality that critical points lie on \new{a finite number of} hyperspheres centered at data points. We further show that critical points are contained in affine linear spaces of specified dimensions. In addition, we observe that the infinite ED discriminant for an affine cone and its dual coincide. 

Recall that $C_X(u)$ denotes the set of critical points for $u$, and that $\mathcal E_X(u)$ is the set of $x$ such that $(x,u)$ lies in the ED correspondence $\mathcal E_X$.

\begin{theorem}[Structure of critical points]\label{thm: structure} Let $X\subseteq \CC^n$ be an irreducible variety. For any $u\in \CC^n$, $\mathcal E_X(u)$ is contained in a finite union of hyperspheres centered at $u$. 
\end{theorem}

\begin{proof} Take an irreducible component $\Lambda$ of $ \mathcal E_X(u) $. There is nothing to prove if $\Lambda$ is a point, so let $\dim \Lambda\ge 1$.

First, assume that $\Lambda$ intersects the smooth locus of $X$. Then a generic point $x\in \Lambda$ is smooth both in $\Lambda$ and in $X$. There is a smooth local \new{parametrization} $x(s)$ of $\Lambda$ in a neighborhood around $x(0)=x$, where $s\in \CC^{\dim \Lambda}$. We have either that $\new{d_u(x(s))}=0$ for every $s$, implying that $x(s)$ lies in the hypersphere of radius $0$ around $u$, or the parametrization takes the form 
\begin{align}
    x(s)=u+r(s)\varphi(s),
\end{align}
where $r(s)\in \CC$ is nonzero around \new{$s=0$ and} $\varphi(s)$ takes values in the unit hypersphere. We show that around $s=0$, $r(s)$ is a constant.

Since $x(s)$ are critical points for $u$, we have $x(s)-u\in N_{x(s)} X$. By assumption of smoothness, $\frac{\partial }{\partial s_i}x(s)$ span the tangent plane of $\Lambda$ at $x(s)$, and therefore  $x(s)-u\in N_{x(s)} X$ is equivalent to 
\begin{align}\label{eq: x x prim}
   \Big(\frac{\partial}{\partial s_i}x(s)\Big) \cdot  (x(s)-u)=0 \quad \textnormal{for each }\new{s_i}.
\end{align}
By the chain rule, \eqref{eq: x x prim} is equivalent to 
\begin{align}\label{eq: r varphi}
   r(s) \Big(\frac{\partial }{\partial s_i}r(s)\Big)\varphi(s)\cdot \varphi(s)+ r(s)^2\Big(\frac{\partial}{\partial s_i}\varphi(s)\Big)\cdot \varphi(s)=0\quad \textnormal{for each }\new{s_i}.
\end{align}
We may factor out $r(s)$, since it is \new{nonzero}. By construction, we have $\varphi(s)\cdot \varphi(s)=1$, which implies $\big(\frac{\partial}{\partial s_i} \varphi(s) \big)\cdot \varphi(s) =0$. Plugging these equalities into \eqref{eq: r varphi}, \new{we get} $\frac{\partial }{\partial s_i}r(s)=0$ for each \new{$s_i$}. 

Second, assume that $\Lambda$ is entirely contained in the singular locus of $X$. Consider a Whitney stratification of $X$ defined by the filtration $F_i$. A Zariski dense open subset of $\Lambda$ is contained in $X_k=F_k\setminus F_{k-1}$ for some $k$. \new{A} generic point $x\in \Lambda$ is smooth in both $\Lambda$ and $X_k$. If we can prove that $x-u\in N_x \:\new{X_k}$, then the above deduction applies and shows that $\Lambda$ lies in a hypersphere around $u$.

The ED correspondence is the Euclidean closure of $E_X$, the set of pairs $(x,u)$ such that $x\in X$ is smooth and is a critical point of $u\in \CC^n$. Therefore, given $x\in \mathcal E_X(u)$, there is a sequence $(x_n,u_n)\in E_X$ converging in the Euclidean topology to $(x,u)$. Whitney's condition A says that $\lim T_{x_n}\:X$ contains $T_x\:X_k$ in the limit. This means that $N_x\:\new{X_k}$ contains the limit of $N_{x_n}\:X$. Then $x_n-u_n\in N_{x_n}\:X$ implies that, in the limit, we have $x-u\in N_x\: X_k$.  
\end{proof}

\new{The result does not}, however, imply that the projection from the offset hypersurface
\begin{align}\begin{aligned}\label{eq: proj offset}
    \mathcal O_X &\to \CC^n,\\
    (u,\epsilon)& \mapsto u
\end{aligned}
\end{align}
is finite-to-one. \new{There may exists data points $u$ such that the ED polynomial $f_u(\epsilon)$ is identically zero in $\epsilon$, even though $\mathcal E_X(u)$ is contained in finitely many hyperspheres.} This is demonstrated in the next example, which was suggested to us by Luca Sodomaco. \new{At first glance, this may appear to contradict \Cref{thm: structure}. Nevertheless, we explain why this is not the case.} We remark that this example considers a different squared distance function $d_u(x)$, \new{although} the principle of \Cref{thm: structure} remains unchanged.

\begin{example}[Cone over twisted cubic]\label{ex: twisted} Let $X$ be the affine cone over twisted cubic in $\CC^4$. It is cut out by 
\begin{align}    \mathrm{rank}\begin{bmatrix} x_0 & x_1 & x_2 \\
    x_1 & x_2 & x_3
    \end{bmatrix}\le 1.
\end{align}
With respect to the Bombieri-Weyl squared distance
\begin{align}
    d_u^{\mathrm{BM}}(x): = (x_0-u_0)^2 + 3(x_1-u_1)^2 + 3(x_2-u_2)^2 + (x_3-u_3)^2, 
\end{align}
the ED polynomial $f_u(\epsilon)$ \new{(defined with respect to $d_{x}^{\mathrm{BM}}$)} is identically zero \new{in $\epsilon$ for a degree-2 surface $Z$ of data points. Elements of $Z$} have one-dimensional fibers over \eqref{eq: proj offset}. Take $u\in Z$ and $(u,\epsilon)\in \mathcal{O}_{X}$ (defined with respect to $d_x^{\mathrm{BM}}$) for any fixed $\epsilon\in \CC$. By Chevalley's theorem, there is a sequence \new{$(x_n,u_n)$ in $\mathcal{E}_{X}$ (defined with respect to $d_{x}^{\mathrm{BM}}$) such that $u_n\to u$ and $d_{x_n}^{\mathrm{BM}}(u_n)=\epsilon_n^2\to \epsilon^2$.} However, there is no guarantee that $x_n$ converges. Indeed, \new{\Cref{thm: structure} implies} that $x_n$ goes to infinity\new{, and that for generic $\epsilon$, there is no critical point $x\in \CC^n$ at complex distance $\epsilon$ to $u$.}
\hfill$\new{\rotsymbol\,}$\end{example}

Under the assumption that $X$ meets the isotropic quadric $Q_\infty$ transversally, there is an alternative proof for \Cref{thm: structure}. By \cite[Corollary 4.3.7]{sodomaco2020distance}, the top coefficient of the ED polynomial is constant \new{in $u$}. Therefore, the projection map \eqref{eq: proj offset} is finite-to-1, and to each data point there are finitely many associated squared complex distances in the ED correspondence.

\begin{corollary}\label{cor: center hyp} Let $X\subseteq \CC^n$ be an irreducible variety. Then $\overline{C_X(c)}=X$ if and only if $c$ is the center of a hypersphere $S$ containing $X$.  
\end{corollary}

\begin{proof} $ $

$\Leftarrow)$ We can locally parametrize $X$ as $x(s)=c+\new{q}(s)$ for $s\in \CC^{\dim X}$ around a smooth point. Then $\new{q}(s)\cdot \new{q}(s)$ is constant, and we have $\big(\frac{\partial}{\partial s_i}x(s)\big)\cdot \new{q}(s)=\big(\frac{\partial}{\partial s_i}\new{q}(s)\big)\cdot \new{q}(s)= 0$. In particular, $\new{q}(s)\in N_{x(s)}X$, showing that $x(s)-c\in N_{x(s)}X$ for each $s$. 

$\Rightarrow)$ By \Cref{thm: structure} and the irreducibility of $X$, $\overline{C_X(c)}$ is contained in a hypersphere around $c$. 
\end{proof}

Critical points exhibit more structure than merely lying on hyperspheres centered at the data point: they are contained in affine linear spaces with specified normal spaces and dimensions as follows.

\begin{lemma}\label{le: diff eq} Let \new{$C,X\subseteq \CC^n$ be irreducible varieties such that $C\subseteq \Sigma_X^\infty$}. For generic $u\in C$, any irreducible component of $\mathcal E_X(u)$ is contained in an affine linear space, whose normal space is $T_u\:C$. 
\end{lemma}

\begin{proof} Let $\Lambda$ be an irreducible component of $\mathcal E_X(u)$. Consider the correspondence
\begin{align}\label{eq: EXC}
    \mathcal E_X^C:=\mathcal E_X\cap \new{( X\times C)}.
\end{align}
We write $\mathrm{pr}_i^{C}$ for the projection onto the $i$-th factor for $i=1,2$. \new{Since $u$ is generic, $\Lambda\times\{u\}$ lies in an irreducible component $Z$ of $\mathcal E_X^C$ such that $\overline{\mathrm{pr}_2^C(Z)}=C$.} By \Cref{le: gen fib}, a generic point $x\in \Lambda$ has that $(x,u)$ is generic in $Z$. Let $W$ be the intersection of a small open (Euclidean) neighborhood (in $\CC^n\times \CC^n$) of \new{$(x,u)$} with $Z$. The projection $B = \mathrm{pr}_2(W)$ is the intersection of an open neighborhood containing $u$ with $C$.

Since $(x,u)$ is generic in $Z$, the differential of $\mathrm{pr}_2^C$ restricted to $Z$ is full rank around this point. By Ehresmann's fibration theorem (\Cref{ss: local par}), we can locally parametrize $W$ 
\begin{align}
    (f(s,t),u(s))\in Z \subseteq \mathcal E_X,
\end{align}
with smooth function $f,u$, where $u(s)$ is a local smooth parametrization of $C$ around $u(0)=u$. We proceed by setting $q(s,t)=f(s,t)-\new{u(s)}$. It follows from \Cref{thm: structure}, that
\begin{align}\label{eq: q}
    q(s,t)\cdot q(s,t) =R(s)
\end{align}
for some \new{smooth} function $R(s)$.

Consider a Whitney stratification of $X$ given by the filtration $F_i$. A dense open subset of the irreducible variety $Y:=\overline{\mathrm{pr}_1^C(Z)}$ is contained in $X_k$ for some $k$. By the argument at the end of the proof of \Cref{thm: structure}, since $u(s)+q(s,t)\in \mathcal E_X(u(s))$, the normal space $N_{f(s,t)}\:Y$ must contain $q(s,t)$. The tangent space $T_{x}\:Y$ is spanned by the partial derivatives of $u(s)+q(s,t)$, implying 
\begin{align}\label{eq: affine hype}
    \Big(\frac{\partial}{\partial s_i}(u(s)+q(s,t))\Big)\cdot q(s,t)=0 \quad \textnormal{ and }\quad  \Big(\frac{\partial}{\partial t_i}(u(s)+q(s,t))\Big)\cdot q(s,t)=0
\end{align}
for each $s_i$ and $t_i$. The latter expression is trivially zero. The former expression becomes
\begin{align}\label{eq: hype}
    \frac{1}{2}\frac{\partial}{\partial s_i}R(s) + \Big(\frac{\partial}{\partial s_i} u(s)\Big)\cdot (f(s,t)-u(s))=0.
\end{align}
The partial derivatives of $u(s)$ span the tangent space $T_{u(s)}\: C$. As a consequence, for fixed $s$, \eqref{eq: hype} defines an affine linear space for $x=f(s,t)$, whose normal is $T_{u(s)} \: C$.
\end{proof}

Finally, we observe that the infinite ED discriminant of an affine cone and its dual coincide:

\begin{proposition}
If $X$ is an affine cone, then 
\begin{align}
    \Sigma_X^{\infty}=\Sigma_{X^*}^{\infty},
\end{align}
where $X^*$ is the dual variety to $X$.
\end{proposition}
\begin{proof}
The statement follows from the duality theorem \cite{draisma2016euclidean}[Theorem $5.2$]. Take a point $u\in\new{\Sigma_X^{\infty}}$. \new{The critical points $x\in C_X(u)$ map to the critical points in $C_{X^{*}}(u)$ via $x\mapsto u-x$.} This extends in a straightforward way to $\mathcal E_X(u)$ and $\mathcal E_{X^*}(u)$.
\end{proof}
%%%%%%%%%%%%%%%%%%%%%%%%%%%%%%%%%%%%%%%%%%%%%%%%%%%%%%%%%%%%%%%%%%%%%%%%%%%%%%%%%%%%%%%%%%%%%%%%%%%%%%%%%%%%%%%%%%%%%%%%%%%%%%%%

\section{\texorpdfstring{Classification of $\Sigma_X^{\infty}$ for curves}{}}\label{s: curves}

Building on \Cref{thm: structure} (and the discussion under it), we derive a complete characterization of infinite ED discriminants for curves $X$ as follows.

\begin{theorem}[Classification for curves]\label{thm: curve class} Let $X$ be an irreducible curve in $\CC^n$.
\begin{enumerate}[label=(\Roman*)]
    \item If $X$ is not contained in a hypersphere, then $\Sigma_X^\infty=\emptyset$.
    \item Otherwise, let $c$ be the center of such a hypersphere. If $V$ is the smallest affine linear space containing \new{$X$ and} $N\:V$ is its normal space, then 
    \begin{align}
     \Sigma_X^\infty=\{c+y: y\in N\: V\}.  
    \end{align}
\end{enumerate}
\end{theorem}

In either case, the infinite ED discriminant of a curve is an affine linear space. This result agrees with \Cref{le: diff eq}, since $T\: \Sigma_X^\infty$ coincides with $N\: V$.

\begin{proof} $ $

$(I).$ We have $c\in \Sigma_X^\infty$ if and only if $\overline{C_X(c)}=X$, since $X$ is one-dimensional. By \Cref{cor: center hyp}, such a $c$ must be the center of a hypersphere containing $X$. By assumption, there are no such hyperspheres.

$(II).$ Let $x(s)$ be a smooth local parametrization of $X$. Again, we use \Cref{cor: center hyp} to see that $c\in  \Sigma_X^\infty$ if and only if $X$ is contained in a hypersphere around $c$. Then
\begin{align}
   \new{d_c(x(s))}= R,
\end{align}
for every $s$ and some $R\in \CC$. After differentiating both sides with respect to $s$, this is equivalent to $\frac{\partial}{\partial s}x(s)\cdot(x(s)-c)=0$ \new{(for every $s$)}. In particular, $\Sigma_X^\infty$ is precisely the set of points $u$ such that $\frac{\partial}{\partial s}x(s)\cdot(x(s)-u)=0$. \new{Note that $u\in \Sigma_X^\infty$ is equivalent to $\frac{\partial}{\partial s}x(s)\cdot (c-u)=0$. We are done if we can show that $T\: V=\mathrm{span}_{s}\{\frac{\partial}{\partial s}x(s)\}$. Observe that $\frac{\partial}{\partial s}x(s)\in T\: V$}, because $x(s)\in V$. Conversely, assume that $\mathrm{span}_{s}\{\frac{\partial}{\partial s}x(s)\}\subsetneq \new{T\:V}$. Then there would exist \new{a} $\new{v}\not\in N\:V$ such that $\frac{\partial}{\partial s}x(s)\cdot \new{v}=0$. As a consequence, there would exist \new{a} $\new{w}$ such that $x(s)\cdot \new{v}+\new{w}=0$. Since $\new{v}\not\in N\: V$, this contradicts the fact that $V$ was the smallest affine linear space containing $X$. It follows that $u= c + y$ for some $y\in N\: V$. 
\end{proof}

\begin{corollary} Let $X$ be an irreducible plane curve in $\CC^2$. 
\begin{enumerate}
    \item $\Sigma_X^\infty$ is a one-dimensional if and only if $X$ is a line of direction $(i,-1)^T$ or $(i,1)^T$. In this case, $\Sigma_X^\infty=X$.
    \item $\Sigma_X^\infty$ is zero-dimensional if and only if $X$ is a circle of nonzero radius, i.e.,
\begin{align}
    X=\{(x_1,x_2):(x_1-c_1)^2+(x_2-c_2)^2=r^2\},
\end{align}
for some center $c=(c_1,c_2)$ and radius $r\neq 0$. In this case, $\Sigma_X^\infty=\{c\}$.
\item $\Sigma_X^\infty$ is empty for all other curves $X$. 
\end{enumerate}
\end{corollary}

%%%%%%%%%%%%%%%%%%%%%%%%%%%%%%%%%%%%%%%%%%%%%%%%%%%%%%%%%%%%%%%%%%%%%%%%%%%%%%%%%%%%%%%%%%%%%%%%%%%%%%

\section{Skew-tubes}\label{s: skew} 
In this section, we define skew-tubes: \new{varieties made up by critical points to data points on a smaller variety.} We characterize skew-tubes in terms of differential equations and then specialize to $\CC^3$, where we provide an algorithm for computing skew-tubes around curves. The algorithm is used to compute several examples. Finally, we discuss optimization properties, focusing on ED degrees of skew-tubes. We demonstrate that these varieties have low ED degrees in relation to their degrees and that their structure can be exploited to accelerate the computation of critical points greatly. 

\begin{definition}\label{def: skew} An irreducible variety $X\subseteq \CC^n$ is a \textit{skew-tube} around $C$ if the following holds: $C$ is an irreducible variety of smaller dimension than $X$ such that $C_X(u)$, for generic $u\in C$, are ($\dim X-\dim C$)-dimensional and their union is dense in $X$. 
\end{definition}

Any skew-tube $X$ around $C$ has that $C\subseteq \Sigma_X^\infty$. Our first result on skew-tubes is a convenient formulation in terms of differential equations for local \new{parametrizations}. 

\begin{lemma}\label{le: skew} Let $C,X\subseteq \CC^n$ be irreducible varieties such that $\dim C<\dim X$. Then $X$ is a skew-tube around $C$ if and only if the following conditions hold: (1) there is a local smooth parametrization $c(s)$ around \new{any generic} point $c(0)=c\in C$ and a function $q(s,t)$ such that $c(s)+q(s,t)$ is a local smooth parametrization around a smooth point of $X$, where $(s,t)\in \CC^{\dim C}\times \CC^{\dim X-\dim C} $, (2) $q(s,t)\cdot q(s,t)=R(s)$ for some smooth $R$, and (3)
\begin{align}\label{eq: cond1}
    \Big(\frac{\partial}{\partial s_i}  (c(s),R(s))\Big)\cdot (q(s,t),\frac{1}{2})=0 \quad \textnormal{for every }s_i.
\end{align}

Further, if these condition are satisfied, then $c(s)+q(s,t)\in C_X(c(s))$ for every small enough $s$ and $t$.  
\end{lemma}

\begin{proof} This result essentially follows directly from the proof of \Cref{le: diff eq}. 

$\Leftarrow)$ The point $c(s)+q(s,t)$ is critical \new{to} $c(s)$ if \new{$q(s,t)$} is (algebraically) orthogonal to the partial derivatives of $c(s)+q(s,t)$, which span the tangent spaces of $X$. \new{The orthogonality} follows from assumptions (2) and (3). 

$\Rightarrow)$ The condition \eqref{eq: cond1} is \eqref{eq: hype}. Indeed, $C$ is \new{a subvariety of $\Sigma_X^\infty$, and as in the proof of \Cref{le: diff eq},} we get the parametrization $c(s)+q(s,t)$ of $X$ that satisfies the conditions (2) and (3).
\end{proof}

%%%%%%%%%%%%%%%%%%%%%%%%%%%%%%%%%%%%%%%%%%%%%%%%%%%%%%%%%%%%%%%%%%%%%%%%%%%%%%%%%%%%%%%%%%%%%%%%%%%%%%%%%%%%%%%%%%%%%%%%%%%%%%%%%%%%%%%%%%%%%%%%%%%%%%%%%%%%

\subsection{\texorpdfstring{Surfaces in $\CC^3$}{}}\label{ss: surfaces R3} 
\begin{figure}
    \centering
    \includegraphics[width=0.4\textwidth]{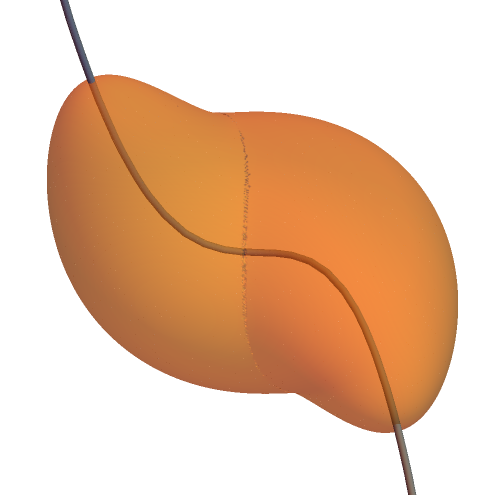}\hspace{2em}
   \includegraphics[width=0.4\textwidth]{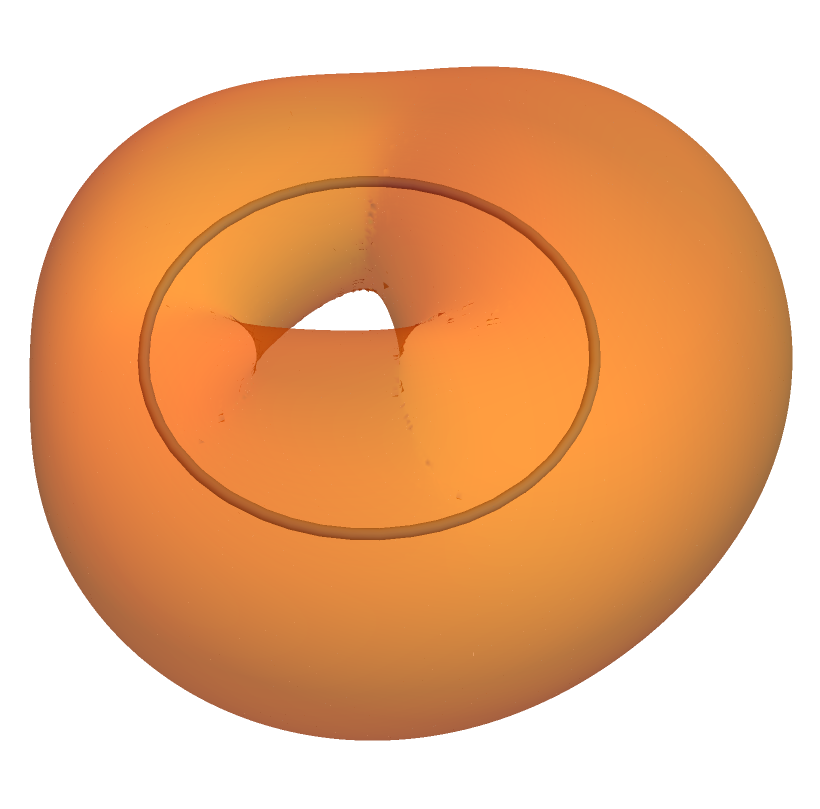}
    \caption{Surface\new{s (orange)} whose infinite ED discriminants contain curves (gray). On the left: the degree-10 surface (T2) and a twisted cubic (\Cref{ex: twisted skewtubes}). On the right: the degree-12 surface (C2) and a circle (\Cref{ex: circle skewtubes}).}
    \label{fig: skewtubes}
\end{figure}

We characterize all surfaces in $\CC^3$ that are skew-tubes around curves via explicit constructions. In \Cref{s: examples}, we encountered such surfaces in the form of surfaces of revolution and tubes. Here, \new{we establish other kinds of surfaces,} including the ones shown in \Cref{fig: skewtubes}. 

Consider an irreducible curve $\mathcal U\subseteq \CC^3 \times \CC$ with coordinates $(c,R)$. Locally, $\mathcal U$ is cut out by 3 polynomials whose Jacobian is of full rank at generic $(c,R)\in \mathcal U$. The tangent space $T_{(c,R)} \:\mathcal U$ is the kernel of the Jacobian, spanned by a single vector defined via Cramer's rule: its coordinates are determinants of submatrices of the Jacobian. This is a choice of \textit{local tangent} $\new{T}_{(c,R)}^{\mathrm{loc}} \:\mathcal U$. Also, let $\mathrm{pr}_1:\mathcal U\to \CC^3$ denote the projection onto the first factor.

\begin{theorem}[Construction of skew-tubes]\label{thm: surf class Imp} Let $\mathcal U\subseteq \CC^3\times \CC$ be an irreducible curve that is cut out locally by $f_1,f_2,f_3\in \CC[c_1,c_2,c_3,R]$\new{, which} determine a local tangent $T_{(c,R)}^{\mathrm{loc}}\:\mathcal U$. Assume that $C=\overline{\mathrm{pr}_1(\mathcal U)}$ is a curve. Define the ideal $I^{\mathcal U}\subseteq S:=\CC[x_1,x_2,x_3,c_1,c_2,c_3,R]$ as the sum of the defining ideal of $\mathcal U$ and the ideal generated by the conditions
\begin{align}\label{eq: extra gens}
   (x-c)\cdot (x-c)=R, \quad T_{(c,R)}^{\mathrm{loc}}\:\mathcal U\cdot \new{\Big(x-c,\frac{1}{2}\Big)}=0.
\end{align}
Let $I_{\mathrm{sat}}^{\mathcal U}$ denote the saturation of $I^{\mathcal U}$ with respect to the ideal generated by $T_{(c,R)}^{\mathrm{loc}}\:\mathcal U=0$. Then any component $X$ of the zero locus of
\begin{align}
   I_{\mathrm{elim}}^{\mathcal U}:=I_{\mathrm{sat}}^{\mathcal U}\cap \CC[x_1,x_2,x_3]
\end{align}
satisfies $C\subseteq\Sigma_X^\infty$. 

Moreover, for any surface $X$ in $\CC^3$ that is a skew-tube around a curve $C$, there is an irreducible curve $\mathcal U\subseteq \CC^4$ with $C=\overline{\mathrm{pr}_1(\mathcal U)}$ such that $X$ is a component of the zero locus of $I_{\mathrm{elim}}^{\mathcal U}$.
\end{theorem}

\begin{example} The algorithm does not always produce a two-dimensional variety\new{,} and $\overline{\mathrm{pr}_1(\mathcal U)}$ is not always an irreducible component of $\Sigma_X^\infty$. To see this, consider $\mathcal U=\mathcal V(\langle c_1+ic_2,c_3,R\rangle)$. Its tangent \new{at any point} is $(i,-1,0,0)$. The ideal $I_{\mathrm{elim}}^{\mathcal U}$ is generated by $x_3=0,ix_1-x_2=0$, and $\Sigma_X^\infty$ is a plane by \Cref{thm: curve class}.  
\hfill$\new{\rotsymbol\,}$\end{example}

\begin{example}\label{ex: isotropic as skew} Let $\mathcal U$ be the curve $c_1^2+c_2^2+1=0,c_3=1,R=0$. For each $(c,0)\in \mathcal U$, the zero locus in $x$ of \eqref{eq: extra gens} is the union of two \new{lines} parametrized by $\lambda$:
\begin{align}
    \lambda( c_1, c_2, \pm 1).
\end{align}
Taking the union as $c$ varies, we get the isotropic quadric $Q$. Indeed, the isotropic quadric is a skew-tube around infinitely many different curves inside $Q$.
\hfill$\new{\rotsymbol\,}$\end{example}

 In the proof, $\mathrm{pr}_{i_1,\ldots,i_j}$ is the projection from a variety in $A_1\times \cdots \times A_m$ onto $A_{i_1}\times \cdots \times A_{i_j}$. Also, $(T_{(c,R)}^{\mathrm{loc}}\:\mathcal U)_{1:3}$ denote the first three coordinates of the local tangent, while $(T_{(c,R)}^{\mathrm{loc}}\:\mathcal U)_4$ is the fourth coordinate.

\begin{proof} Consider \new{first} $(c,R)\in \mathcal U$ with $T_{(c,R)}^{\mathrm{loc}}\:\mathcal U\neq 0$. If $(T_{(c,R)}^{\mathrm{loc}}\:\mathcal U)_{1:3}=0$, then $(T_{(c,R)}^{\mathrm{loc}}\:\mathcal U)_4\neq 0$ and there is no solution in $x$ to \eqref{eq: extra gens}. If $(T_{(c,R)}^{\mathrm{loc}}\:\mathcal U)_{1:3}\neq 0$, then we proceed to show that there is a one-dimensional solution \new{set} in $x$ to \eqref{eq: extra gens}. The constraint on the right defines an affine plane $P$ that is not contained in the isotropic quadric (no plane is). This means that we can parametrize $P$ as $\lambda_1p_1+\lambda_2p_2+p_3$ with $p_1\cdot p_1,p_2\cdot p_2\neq 0$. Plugging this into the condition on the left, there is a one-dimensional solution set in $\lambda_1,\lambda_2$\new{, which suffices. 

There are finitely many $(c,R)$ such that $T_{(c,R)}^{\mathrm{loc}}\:\mathcal U= 0$, and for these} there is a two-dimensional solution \new{set} in $x$ to \eqref{eq: extra gens}. These constitute components of $V(I^{\mathcal U})$ that are removed during saturation: they are not contained in $V(I_{\mathrm{sat}}^{\mathcal U})$. This means that for fixed $(c,R)\in \mathcal U$, the solution set in $x$ to $(x,c,R)\in V(I_{\mathrm{sat}}^{\mathcal U})$ is at most a \new{curve}, and for generic $\new{(c,R)}\in \mathcal U$, the solution set is a \new{degree-2} curve. 

The above shows that $V(I_{\mathrm{sat}}^{\mathcal U})$ is two-dimensional. Let $Y_1,\ldots,Y_k$ be the irreducible components of $V(I_{\mathrm{sat}}^{\mathcal U})$. We have $\mathcal U=\overline{\mathrm{pr}_{2,3}(Y_i)}$, at least for some $i$. By the semi-upper continuity of dimensions~\cite[Chapter 1, Section 8, Corollary 3]{mumford1999red}, \new{we deduce that} $k\le 2$ and $\mathcal U=\overline{\mathrm{pr}_{2,3}(Y_i)}$ for all $i$.

Let $Y=\new{Y_i}$ and define $X$ as the closure of the image of $(x,c,R)\mapsto x$ restricted to $Y$. By Ehresmann's fibration theorem (\Cref{ss: local par}), we can locally parametrize $Y$ via smooth functions $(x(s,t),c(s),R(s))$, with $(s,t)\in \CC^1\times \CC^1$. \new{The parametrization} satisfies 
\begin{align}\label{eq: two cond reform}
    \new{d_{c(s)}(x(s,t))}=R(s), \quad \Big(\frac{\partial}{\partial s}( c(s),R(s))\Big) \cdot (x(s,t)-c(s),\frac{1}{2}),
\end{align}
since $\frac{\partial}{\partial s}( c(s), R(s) )$ is the tangent direction of $\mathcal U$ at $(c(s),R(s))$. The tangent space $\new{T_{x(s,t)}}\: X$ is spanned by the partial derivatives of $x(s,t)$. As in the proof of \Cref{le: diff eq}, $\new{\frac{\partial}{\partial s}x(s,t)}\cdot (x(s,t)-c(s))=0$ is equivalent to the right condition of \eqref{eq: two cond reform} and $\frac{\partial}{\partial t}x(s,t)\cdot (x(s,t)-c(s))=0$ is equivalent to the left condition. This shows that $x(s,t)\in C_X(c(s))$ for \new{small $s,t$}, which suffices to show that $C\subseteq \Sigma_X^\infty$.

For the second part \new{of the statement, let $Z$ be an irreducible component of $\mathcal E_X^C$ (defined in \eqref{eq: EXC}) such that $X=\overline{\mathrm{pr}_1(Z)}$ and $C=\overline{\mathrm{pr}_2(Z)}$. The dimension of $Z$ is 2. Consider the map
\begin{align}
    \pi: Z \to \CC^3\times \CC, \quad (x,c)\mapsto (c,d_c(x)). 
\end{align}
By \Cref{thm: structure}, for each $c\in C$, there are finitely many associated squared distances $d_c(x)$. Then the fiber dimension theorem~\cite[Chapter 1, Section 8, Corollary 1]{mumford1999red} implies that the generic fiber of $\pi$ is one-dimensional. The} closure of the image of $\pi$ is a curve $\mathcal U$. In particular, $C=\overline{\mathrm{pr}_1(\mathcal U)}$.

As in \Cref{le: skew}, there is a local parametrization $x(s,t)=c(s)+q(s,t)$ of $X$ around \new{$x(0,0)=x$ for which \eqref{eq: cond1} holds}. The local tangent of $\mathcal U$ at $(c(s),R(s))\in \mathcal U$ coincides with $\frac{\partial}{\partial s}(c(s),R(s))$, by construction. If $W$ \new{is} the closure of the image of $(x,c)\mapsto (x,c,\new{d_c(x)})$ restricted to $Z$, then $W$ is contained in $V(I_{\mathrm{sat}}^{\mathcal U})$. It is clear that $V(I_{\mathrm{sat}}^{\mathcal U})$ cannot equal $\new{\CC^3\times \mathcal U}$, and therefore $V(I_{\mathrm{elim}}^{\mathcal U})$ is at most two-dimensional. Then $V(I_{\mathrm{elim}}^{\mathcal U})$ contains $X = \overline{\mathrm{pr}_{1}(W)}$ as a component. 
\end{proof}

\begin{example}\label{ex: circle skewtubes} Consider the unit circle $C\subseteq \CC^3$ in the $c_1c_2$-plane. With $f_1(c,R)=c_1^2+c_2^2-1$ and $f_2(c,R)=c_3$, we examine skew-tubes for the following choices of $f_3$: (C1) $R-(c_1c_2-1/4)$ \new{and} (C2) $R-((c_1^3+c_2^3)/4+2)$. We compute them in \texttt{Macaulay2} using \Cref{fig: code skew}. The skew-tubes around \new{$C$} are as follows: (C1) a degree-8 surface defined by a polynomial of 55 terms (\Cref{fig: teaser}, right) and (C2) a degree-12 surface defined by a polynomial of 226 terms (\Cref{fig: skewtubes}, right). 
%https://www.desmos.com/3d/4a6a81b455 %https://www.desmos.com/3d/2eabab648d%https://www.desmos.com/3d/79ec79ad5e
\hfill$\new{\rotsymbol\,}$\end{example}

\begin{example}\label{ex: twisted skewtubes} Consider the twisted cubic $C\subseteq \CC^3$ parametrized by $s\mapsto (s,s^2,s^3)$. We examine the following choices of parametrizations of $R(s)$: (T1) $s^2$ and (T2) $1-s^2$. \new{The skew-tubes around the twisted cubic are as follows:} (T1) a degree-10 surface defined by a polynomial of 97 terms (\Cref{fig: teaser}, left) and (T2) a degree-10 surface defined by a polynomial of 123 terms (\Cref{fig: skewtubes}, left).
%https://www.desmos.com/3d?lang=sv-SE
% R = t % https://www.desmos.com/3d/2beeabc28d
% R = t^2 % https://www.desmos.com/3d/4a9b237172
% R = 1-t^2 % https://www.desmos.com/3d/b22541876c
% R = t^3-6t^2+4t+12 % https://www.desmos.com/3d/0685e26c3f
% R = 2*(x-1/2)*(x+1/2)*(x-1) https://www.desmos.com/3d/a5a765cd41
\hfill$\new{\rotsymbol\,}$\end{example}

\begin{lstlisting}[language=Macaulay2 , caption={\texttt{Macaulay2} code for the computation of skew-tubes.} , label={fig: code skew}]    
S = QQ[x_1,x_2,x_3,c_1,c_2,c_3,R]
I = -- * Ideal of U * -- 
If = -- * Ideal of f_1,f_2,f_3 * --

jacIf = diff(matrix{{c_1,c_2,c_3,R}}, transpose gens If)
tangent = matrix{{det submatrix(jacIf,,{1,2,3}),
                - det submatrix(jacIf,,{0,2,3}),
                det submatrix(jacIf,,{0,1,3}),
                - det submatrix(jacIf,,{0,1,2})}}

C = matrix{{c_1},{c_2},{c_3}}
X = matrix{{x_1},{x_2},{x_3}}
IU = I + ideal(transpose( X - C )*( X - C ) - R,
                tangent*( X - C || matrix{{1/2}}))
IUsat = saturate(IU, ideal(tangent))
IUelim = eliminate(toList(c_1,c_2,c_3,R), IUsat)

IX = decompose IUelim
\end{lstlisting}

\subsection{Computation of critical points}\label{ss: opt}
The ED degree of a surface of revolution $X$ is closely related to the ED degree of the curve $C$ being revolved. Indeed, in any plane through the axis of revolution, the smooth points of the surface of revolution have that their normal lines are inside that plane. Then, for a generic data point $u$, all critical points lie in the plane $P$ spanned by $u$ and the axis of revolution. The intersection $X\cap P$ is the union of two copies of the original planar curve (one reflected along the axis) by construction. It follows that
\begin{align}
    \mathrm{EDdeg}(X)= 2\cdot \mathrm{EDdeg}(C).
\end{align}
A straightforward argument also shows that this relation also hold for tubes $X$ around $C$. 

This naturally raises the question whether the ED degree of a skew-tube around a curve is twice the ED degree of the curve. Unfortunately, this is not the case in general. Instead, we prove an upper bound for their ED degrees. In this direction, we first recall from \cite[Proposition 2.6]{draisma2016euclidean}, that for any surface $X\subseteq \CC^3$, 
\begin{align}\label{eq: upper bound ED}
    \mathrm{EDdeg}(X)\le (\deg X)^3-(\deg X)^2+ \deg X,
\end{align}
with equality if $X$ is generic. 

\new{Our skew-tube specialized version of this upper bound goes as follows.

\begin{proposition}[Optimization for skew-tubes]\label{prop: skew eddeg} Let $\mathcal U\subseteq \CC^3\times \CC$ be an irreducible curve that is cut out locally by $f_1,f_2,f_3\in \CC[c_1,c_2,c_3,R]$. Consider $I_{\mathrm{elim}}^{\mathcal U}$ as in~\Cref{thm: surf class Imp}. For a surface component $X$ of $V(I_{\mathrm{elim}}^{\mathcal U})$, we have  
\begin{align}\label{eq: critical skew}
    \mathrm{EDdeg}(X) \le d_1d_2d_3\big(2(d_1+d_2+d_3)-3\big),
\end{align}
where $\deg f_i=d_i$. 
\end{proposition}}

\new{In order to prove this result, we first establish skew-tube specific critical equations. Recall that $I:J^\infty$ denotes saturation of $I$ with respect to $J$.}

\begin{lemma} Let $\mathcal U\subseteq \CC^3\times \CC$ be an irreducible curve that is cut out locally by $f_1,f_2,f_3\in \CC[c_1,c_2,c_3,R]$. Assume that $R$ is \new{nonconstant} in $\mathcal U$ and consider $I_{\mathrm{elim}}^{\mathcal U}$ as in~\Cref{thm: surf class Imp}. For a data point $u\in \CC^3$, we define
\begin{align}\label{eq: Fu}
    F(u):=\new{\big(}(T_{(c,R)}^{\mathrm{loc}}\mathcal U)_4\new{\big)}^2 \new{d_c(u)} - 4 R \big((T_{(c,R)}^{\mathrm{loc}}\mathcal U)_{1:3}\cdot (u-c)\big)^2.
\end{align}
If $X=V(I_{\mathrm{elim}}^{\mathcal U})$ is an irreducible surface, and $u$ generic, \new{then every critical point $x$ to $u$ is of the form 
\begin{align}\label{eq: lamstar}
   x= c+\lambda^*(u-c),\quad \textnormal{where} \quad \lambda^*:= -\frac{(T_{(c,R)}^{\mathrm{loc}}\mathcal U)_4}{2(T_{(c,R)}^{\mathrm{loc}}\mathcal U)_{1:3}\cdot (u-c)}
\end{align}
for some $(c,R)$ in the zero locus of}
\begin{align}\label{eq: I crit U}
    I_{\mathrm{Crit}}^{\mathcal U}(u):= \langle f_1,f_2,f_3, F(u) \rangle : \Big(\langle (T_{(c,R)}^{\mathrm{loc}}\mathcal U)_{1:3}\cdot (u-c) \rangle \cap \langle(T_{(c,R)}^{\mathrm{loc}}\mathcal U)_4 \rangle \Big)^\infty.
\end{align}

\new{Moreover, in this case, the map $V(I_{\mathrm{Crit}}^{\mathcal U}(u))\to C_X(u), (c,R)\mapsto x = c +\lambda^*(u-c)$ is a bijection.}
\end{lemma}
As follows from the proof of \Cref{thm: surf class Imp}, $V(I_{\mathrm{elim}}^{\mathcal U})$ consists of at most two components. In all of our examples (\Cref{ex: circle skewtubes,ex: twisted skewtubes}), it is an irreducible surface. In the reducible case, the critical points are \new{still of the form \eqref{eq: lamstar} with $(c,R)$ contained in $V(I_{\mathrm{Crit}}^{\mathcal U}(u))$}.

\begin{proof} \new{Since $u$ is generic, any $x\in C_X(u)$ is generic in $X$. By the proof of \Cref{thm: surf class Imp}, for such an $x$, there is a $(c,R)\in \mathcal U$ with $x\in C_X(c)$ and \eqref{eq: extra gens} holds. Both $x-u$ and $x-c$ are contained in $N_x\: X$. In particular, $x$ lies on the line spanned by $c$ and $u$, \new{meaning} $x=c+\lambda(u-c)$ for some $\lambda\in \CC$. Then \eqref{eq: extra gens} becomes}
\begin{align}\label{eq: skew critical}
    \lambda^2(u-c)\cdot (u-c)=R,\quad T_{(c,R)}^{\mathrm{loc}}\:\mathcal U \cdot (\lambda(u-c),\frac{1}{2})=0.
\end{align}
\new{Assuming that $(T_{(c,R)}^{\mathrm{loc}}\mathcal U)_{1:3}\cdot (u-c)\neq 0$, we} can solve for $\lambda$ in the second condition of \eqref{eq: skew critical} to get \new{$\lambda^*$ as in \eqref{eq: lamstar}.}
Plugging \new{$\lambda^*$} into the first condition of \eqref{eq: skew critical}, we get the expression $F(u)=0$. \new{Since $u$ is generic and $R$ is nonconstant in $\mathcal U$, $(T_{(c,R)}^{\mathrm{loc}}\mathcal U)_{1:3}\cdot (u-c)=0$ and $(T_{(c,R)}^{\mathrm{loc}}\mathcal U)_4=0$ have no common solution in $(c,R)\in \mathcal U$.} This means that if $(c,R)$ is such that $(T_{(c,R)}^{\mathrm{loc}}\mathcal U)_{1:3}\cdot (u-c)=0$, the second condition of~\eqref{eq: skew critical} is not satisfied. \new{We have shown that the critical point $x$ equals $c + \lambda^*(u-c)$ for some $(c,R)\in V(I_{\mathrm{Crit}}^{\mathcal U}(u))$.}

\new{Conversely, take some $(c,R)\in V(I_{\mathrm{Crit}}^{\mathcal U}(u))$. For such $(c,R)$, $\lambda^*$ is well-defined and we wish to show that $x = c + \lambda^*(u-c)$ is critical to $u$. Since $x$ solves \eqref{eq: skew critical}, the proof of \Cref{thm: surf class Imp} says that $x\in \mathcal E_X(u)$. We need to prove that $x$ is smooth in $X$. First, we observe that no $(c,R)\in \mathcal U$ solves the critical equations for all $u$, which is a consequence of the following fact: the expression $v\cdot v/(w\cdot v)^2$ is nonconstant over $v$ for fixed nonzero $w$ (here, $v$ represents $u-c$ and $w$ represents $(T_{(c,R)}^{\mathrm{loc}}\mathcal U)_{1:3} $ in the equation $F(u)=0$). It follows that $(c,R)$ is generic in $\mathcal U$. Consider the cone through $c$ over the points $x$ that satisfy \eqref{eq: extra gens} (generic such $x$ are smooth in $X$). The data point $u$ is generic in this cone, which meets the singular locus of the skew-tube in at most finitely many points. The line through $u$ and $c$, therefore, does not meet the singular locus.}   
\end{proof}

\begin{proof}[\new{Proof of \Cref{prop: skew eddeg}}] \new{First, assume that $R$ is nonconstant in $\mathcal U$. We estimate the cardinality of $I_{\mathrm{Crit}}^{\mathcal U}(u)$ from above}.  By~\cite[Section I.7, Theorem 7.7]{hartshorne2013algebraic}, the zero locus of $I_{\mathrm{Crit}}^{\mathcal U}(u)$ is a finite number of points, bounded above by the product of the degrees of the generators. The degree of $ T_{(c,R)}^{\mathrm{loc}}\:\mathcal U$ is at most $d_1+d_2+d_3-3$ by construction, implying that~$F(u)$ is at most \new{of} degree $2(d_1+d_2+d_3)-3$ and the upper bound follows. 

In the case that $R$ is constant in $\mathcal U$, then $X$ is an $\epsilon$-offset hypersurface, and the ED degree of $X$ is twice that of $C$. In particular, an upper bound for the ED degree is $2d_1d_2d_3$, agreeing with~\eqref{eq: critical skew}.
\end{proof}

The ED polynomials of the surfaces in~\Cref{ex: twisted skewtubes,ex: circle skewtubes} could not be computed in a reasonable time on our system (see hardware details in~\Cref{s: Comp}). Instead, we conjecture their ED degrees by making a choice of a pseudo-random data point and computing the critical equations. We consider the \textit{standard critical ideal}~\cite[Section 2]{draisma2016euclidean} and the \textit{skew-tube critical equations} of \eqref{eq: I crit U}. In the example below, we see that \Cref{prop: skew eddeg} provides a much tighter bound than~\eqref{eq: upper bound ED}. Moreover, by exploiting the skew-tube structure, the computation of critical points is sped up by factors ranging from 300 to \new{250 000}.

\begin{example}\label{ex: skew ED} For the surfaces in \Cref{ex: twisted skewtubes,ex: circle skewtubes}, the upper bounds for the ED degree obtained from \eqref{eq: upper bound ED} are: (C1) 456, (C2) 1 596, (T1) 910, and (T2) 910. The improved upper bounds obtained from \Cref{prop: skew eddeg} are: (C1) 28, (C2) 54, (T1) 28, and (T2) 28. (For the skew-tubes around twisted cubics, we use that the twisted cubic is locally cut out by $c_2-c_1^2=0$, $c_2^2-c_1c_3=0$.) The number of critical points given the pseudo-random data point $u=(5,7,11)$, computed either via \new{the skew-tube critical equations (\Cref{fig: skew EDdeg}) or the standard ones are}: (C1) 10, (C2) 14, (T1) 10, and (T2) 12. The corresponding computation times, skew-tube vs. standard equations, are: (C1) 0.007 \new{seconds} vs.~2 seconds, (C2) 0.009 \new{seconds} vs.~\new{6 minutes}, (T1) 0.006 \new{seconds} vs.~\new{18 minutes}, and (T2) 0.007 \new{seconds} vs.~\new{45 minutes}. 
\hfill$\new{\rotsymbol\,}$
\end{example}

\begin{lstlisting}[language=Macaulay2 , caption={\texttt{Macaulay2} code for the skew-tube critical equations.} , label={fig: skew EDdeg} ]
S = QQ[c_1,c_2,c_3,R]
I = -- * Ideal of U * --
If = -- * Ideal of f_1,f_2,f_3 * --

jacIf = diff(matrix{{c_1,c_2,c_3,R}}, transpose gens If)
tangent = matrix{{det submatrix(jacIf,,{1,2,3}),
                - det submatrix(jacIf,,{0,2,3}),
                det submatrix(jacIf,,{0,1,3}),
                - det submatrix(jacIf,,{0,1,2})}}

C = matrix{{c_1},{c_2},{c_3}}
U = sub(u, S)
Ieq = ideal(tangent_(3)_(0)^2*transpose( U - C )*( U - C ) - 
        4*R*(submatrix(tangent,,{0,1,2})*( U - C ))^2)
Icrit = saturate(I + Ieq,
                intersect(ideal( submatrix(tangent,,{0,1,2}) ),
                ideal(tangent_(3)_(0))))
                
dim Icrit, degree Icrit
\end{lstlisting}

{\small
\bibliographystyle{alpha}
\bibliography{final_revision_main}
}

\end{document}